\documentclass[review,10pt]{elsarticle}
\usepackage[numbers]{natbib}
 \usepackage{mathrsfs}
\usepackage{graphicx}
\usepackage{subfigure}
\usepackage{epstopdf}
\usepackage{amssymb}
\usepackage{color}
\usepackage{amsmath}
\usepackage{float}
\makeatletter\@addtoreset{equation}{section}

\newtheorem{thm}{Theorem}[section]
\newtheorem{cor}[thm]{Corollary}
\newtheorem{lem}[thm]{Lemma}

\newtheorem{rem}[thm]{Remark}

\numberwithin{equation}{section}
\numberwithin{table}{section}
\numberwithin{figure}{section}

\makeatletter\@addtoreset{equation}{section}

\newenvironment{proofed}[1]{\par \textbf{Proof}\quad #1}{\hfill \textbf{} $\Box$ }

\newtheorem{example}{\bf{Example}}[section]

\numberwithin{equation}{section}
\numberwithin{table}{section}
\numberwithin{figure}{section}

\begin{document}
\begin{frontmatter}
\title{Finite difference schemes   for   multi-term  time-fractional mixed  diffusion-wave equations
}

\author[seu]{Zhao-peng Hao}
\ead{haochpeng@126.com}
\author[purdue]{Guang Lin\corref{cor1}}
\cortext[cor1]{Corresponding author. }
\ead{guanglin@purdue.edu}

\address[seu]{Department of Mathematics, Southeast University, Nanjing 210096, P.R.China.}
\address[purdue]{Department of Mathematics $\&$ School of Mechanical Engineering, Purdue University, West Lafayette, IN 47907, USA.}

\begin{abstract}
 The multi-term time-fractional mixed diffusion-wave equations (TFMDWEs) are considered and the numerical method with its  error analysis  is presented in this paper. First, a   $L2$ approximation is proved with  first order accuracy  to the Caputo  fractional derivative of order $\beta \in (1,2).$  Then  the approximation is applied to solve a  one-dimensional TFMDWE   and  an implicit, compact difference scheme is  constructed.   Next, a rigorous error analysis of  the proposed   scheme is carried out by employing the energy method, and it is  proved  to be  convergent with first order accuracy  in time and fourth order in space, respectively.
In addition, some results for   the distributed order  and  two-dimensional extensions  are  also reported in this work. Subsequently, a practical  fast solver   with  linearithmic complexity is   provided with partial diagonalization technique.  Finally, several  numerical examples  are given to demonstrate the accuracy and efficiency of proposed schemes.
\end{abstract}
\end{frontmatter}

\section{Introduction}

In this work, we are concerned with  numerical methods for the multi-term time-fractional mixed diffusion-wave equations (TFMDWEs) in the following form
\begin{eqnarray}\label{multi-dw}
&&\sum_{i=1}^sK_i\  _0^{C}D_t^{\alpha_i}u(\mathbf{x},t)= \Delta u(\mathbf{x},t) +f(\mathbf{x},t),\quad \mathbf{x}\in \Omega,~t\in (0, T],
\end{eqnarray}
where $\Omega$ is spatial domain, $\Delta$ is the Laplace operator,  $0< \alpha_1<\cdots \leq  \alpha_{i}\leq 1<\alpha_{i+1}< \cdots <\alpha_s <2 $ and $K_i>0,$ $i=1,\ldots,s.$ Here $ \  _0^{C}D_t^{\alpha_i}$ denotes the Caputo fractional derivative of order $\alpha_i,$ which will be defined later.

  As a natural extension of the single-term time-fractional partial differential  equations (TFPDEs), e.g. sub-diffusion or diffusion-wave equations, the multi-term TFPDEs are expected to improve the modeling accuracy in depicting the anomalous diffusion process, successfully  capturing power-law frequency dependence\cite{Kelly2008}, adequately  modeling  various types of viscoelastic damping \cite{Chen-Liu2012}. For instance, in \cite{Schumer-Benson2003},  a two-term  mobile and immobile fractional-order diffusion model was proposed to model  the total concentration in solute transport. The kinetic equation with two fractional derivatives of different orders appears
 quite natural when describing sub-diffusive motion in velocity fields \cite{Metzler1998}. The two-term time-fractional telegraph equations  \cite{Orsingher2004}, which can also be regarded as special cases of  Equation \eqref{multi-dw},  govern
 the iterated Brownian motion and the telegraph processes with Brownian
time.

 There are a few mathematical theories  on multi-term TFPDEs.
 Based on an appropriate maximum principle, Luchko \cite{Luchko-2011} derived some priori
estimates for the solution and   established its uniqueness using the Fourier
method and  the method of separation of variables. Daftardar-Gejji and Bhalekar \cite{Daftardar-Gejji2008} considered a multi-term
time-fractional diffusion-wave equation (TFDWE) along with the homogeneous/nonhomogeneous boundary
conditions and have  solved this equation  using the method of separation of variables. Based
on the orthogonal polynomials of the Laguerre type, Stojanovic \cite{Stojanovic2011} found solutions for the
diffusion-wave equation in one dimension with $n$-term time-fractional derivatives, whose orders
belong to the intervals (0, 1), (1, 2) and (0, 2), respectively. By using Luchko's Theorem and the
equivalent relationship between the Laplacian operator and the Riesz fractional derivative, Jiang
et al. \cite{Jiang2012} derived the analytical solutions for the multi-term time-space fractional advection-diffusion
equations. Subsequently, Ding et al. \cite{Ding2013} presented the analytical solutions
for the multi-term time-space fractional advection-diffusion equations with mixed boundary
conditions. Recently, Stojanovic has  analyzed regularity,  existence and uniqueness of the equation with nonlinear source  in  \cite{Stojanovic2013}.

Numerical approximations of multi-term TFPDEs  have also been discussed by some authors.  For the multi-term TFPDEs  whose orders, $\alpha_i$, belong to $(0,1),$ very recently, Jin et al. \cite{Jin2015} developed a fully discrete scheme based  on the standard Galerkin
finite element method in space  combining with a finite difference discretization of the time-fractional derivatives.   For the orders, $\alpha_i$, lying in  $(1,2),$ Chen. et al.  \cite{Chen-Liu2012} presented a finite difference scheme and gave its analysis, following the idea of the method of order reduction  proposed by Sun and Wu \cite{Sun2006}. In addition, an extension of multi-term TFPDEs, distributed order TFPDEs have also been considered (see e.g.  \cite{Ye2014,Morgado2014}).

 However, to the authors' best knowledge, published works on numerical solution of the multi-term TFMDWEs with the fractional orders, $\alpha_i$, lying in (0,2)   are still  limited in spite
of rich literatures on its single-term version \cite{Sanz-Serna1988,Marcos1990,Tang1993,Sun2006,Lin-Xu2007,Li-Xu2009,Alikhanov2015,Wang2015} and ordinary fractional differential equations ( see e.g. \cite{Katsikadelis2009,Yang2010,Doha2011}).    In \cite{Liu2013},  the authors proposed a numerical method by the method of order  reduction   for more general problems, i.e., the fractional orders, $\alpha_i$, belonging to $(0,n)(n>2).$  However,  the corresponding numerical analysis was not provided.
Due to the definition of the
fractional order derivative, which is a nonlocal operator,  the fractional
order derivative requires a longer memory of the solution. The global nature of the non-integer order
derivative makes the design of stable and  accurate   methods more difficult.  It is
a great challenge for memory and storage requirement when all the past history of the solution has to be saved in order
to compute the solution at the current time. Introducing new variables and transforming the original problem into an equivalent system will make the design of accurate and robust  methods  easier. However,  it will introduce   extra computational cost and data storage meanwhile. Naturally, developing high-order
 numerical methods for time-fractional partial differential equations  is an effective approach to overcome this challenge. But the constraint on the regularity of the underlying problem makes  pursing high-order methods in temporal direction very difficult, since such problem may lack high regularity  due to  the existence of singularity of the fractional derivatives near the initial time.

  The main objective of this paper is to develop a numerical algorithm for multi-term  TFMDWEs  and give  its rigorous  numerical analysis.   Numerous works focus on   time fractional partial differential  equations with  the orders either  lying in (0,1) or in (1,2).  To date,  we are not aware of any published papers investigating the   multi-term  TFMDWEs or their generalization, i.e.,  distributed-order equations with the order belonging to  $(0,2).$  We shall give a complete analysis of our scheme.   By the energy method, a priori estimate for the solution of the Dirichlet boundary value problem of the time-fractional mixed diffusion-wave equations has been established. Instead of  using the method of order reduction  in \cite{Liu2013}, we will utilize the  $L2$ approximation, proposed by Oldham and Spanier \cite{Oldham1974}, to discretize  the Caputo derivative of order $\beta \in (0,2)$ directly, and modify it to be suitable for solving an initial value problem. Combining a modified $L2$ approximation with the well-known $L1$ approximation to discretize  time-fractional derivatives, we construct a compact finite difference scheme by applying the second-order central difference quotient to approximate the weighted average of the second-order derivatives in spatial derivation. By the discrete analogue of energy  or fractional Sobolev inequalities,    we provide  error  analysis of our method rigorously, which can be seen as one of  main contributions of this work.    In addition,   generalization of the distributed order equations and   high-dimensional cases   are  also reported in this paper.  To our knowledge, this is the first paper to investigate the distribution order equations with the orders from 0 to 2. This can be seen another main contributions of this paper.  Other contributions also include that  we present  a fast solver for the fully discretized scheme.   By partial diagonalization technique, the resulting matrix equation  is   reduced to the independent   linear systems  with toeplitz-like structure, being easy and convenient  to design fast algorithms.

The rest of this paper is organized as follows.  In Section 2,   a first-order   $L2$ approximation to the  Caputo fractional derivatives is proved  and   the compact difference scheme is derived. The error analysis  of the proposed scheme is showed in  Section 3. Extension to distributed-order version and to the two-dimensional counterpart are also considered in Section 4 and Section 5, respectively.  In Section 6, a fast solver is suggested with partial diagonalizaion technique and divide and conquer strategy. To show the effectiveness of the proposed algorithm, the numerical experiments are performed  to verify the theoretical results in Section 7.    Finally, some remarks and discussions   are offered in the conclusion part to   close  the paper.

\section{Numerical method for the one-dimensional two-term time-fractional  mixed diffusion-wave equation}
\setcounter{equation}{0}

For the sake of simplicity but without loss of  generality,  we  consider the one-dimensional two-term TFMDWEs as follows:
\begin{eqnarray}
&&K_1\  _0^{C}D_t^{\alpha}u(x,t)+K_2\  _0^{C}D_t^{\beta}u(x,t)= \partial_x^{2}u(x,t) +f(x,t),\quad x\in \Omega,~0<t\leq T, \label{problem-1}\\
&&u(x,0)=\phi_0(x),\quad u_t(x,0)=\phi_1(x),\quad x\in \bar{\Omega},\label{problem-2} \\
 &&u(0,t)=\varphi_0(t),\quad u(L,t)=\varphi_1(t),\quad 0< t\leq T ,\label{problem-3}
\end{eqnarray}
where $\Omega=(0,L),$ $K_1,K_2> 0,$   and  $0< \alpha< 1 <\beta< 2.$

In the following analysis of the proposed  numerical method, we assume the problem \eqref{problem-1}-\eqref{problem-3} has a unique and sufficiently smooth solution.
Here it should be pointed out  that Alikhanov   showed stability of solution by priori estimate for the sub-diffusion and diffusion-wave equations in \cite{Alikhanov2010}.  Following his idea,  one can prove the stability  for the solution of the problem \eqref{problem-1}-\eqref{problem-3} by the energy method.
\subsection{Notations and lemmas}

Take an integer $N.$ Let $\Omega_\tau\equiv \{ t_n~|~0\leq n\leq N\}$ be a uniform mesh of the interval $[0,T],$ where $t_n=n\tau,~0\leq n\leq N$ with $\tau=T/N.$   Two lemmas below are needed. The first one  is well-known  $L1$ approximation.
 \begin{lem}(See \cite{Sun2006})\label{dis-lem-sun-1}
For $\alpha \in (0,1), $ suppose $v(t) \in C^2[0,t_n]$ and denote $v^n=v(t_n),$
\begin{equation}
a_0^{(\alpha)}=\frac{1}{\Gamma{(2-\alpha)}}, \quad  a_n^{(\alpha)}=
\frac{1}{\Gamma{(2-\alpha)}}[(n+1)^{1-\alpha}-n^{1-\alpha}],\qquad n\geq 1, \label{coefficients}\end{equation}
 and define the backward difference quotient operator
\begin{equation}\label{eq:def-al}
\delta_t^{\alpha}v^n= \frac{1}{ \tau^{\alpha} }\biggl[a_0^{(\alpha)}v^n- \sum_{k=1}^{n-1} (a_{n-k-1}^{(\alpha)}-a_{n-k}^{(\alpha)} )v^k -a_{n-1}^{(\alpha)}v^0\biggl],\quad n \geq 1.
\end{equation}
It holds that
\begin{equation}
_0^{C}D_t^{\alpha}v(t_n)=\delta_t^{\alpha}v^n +R_1[v(t_n)],
\end{equation}
where
$$|R_1(v(t_n))|\leq c_{\alpha}\max_{0\leq t\leq t_n}|v^{\prime\prime}(t)|\tau^{2-\alpha}. $$
\end{lem}

 In above Lemma $c_{\alpha}$ is a bounded constant, which is  dependent on $\alpha$ but independent of step size  $\tau.$
For the approximation of derivative of fractional order $\beta \in(1,2),$ we have the following result.
\begin{lem}\label{dis-lemma-2}
For $\beta \in(1,2), $  suppose  $v(t)\in C^{2}[0,t_n]\cap C^{3}(0,t_n], $ and $|v^{\prime\prime\prime}(t)|\in L_1(0,t_n].$  Denote $v^n=v(t_n),$
\begin{equation}
b_0^{(\beta)}:=a_{0}^{(\beta-1)}=\frac{1}{\Gamma{(3-\beta)}},\quad b_n^{(\beta)}:=a_n^{(\beta-1)}=
\frac{1}{\Gamma{(3-\beta)}}[(n+1)^{2-\beta}-n^{2-\beta}],\quad n\geq 1 \label{coefficients-2}\end{equation}
and
\begin{eqnarray}\label{eq:def-be}
&&\Delta_t^{\beta}v^n=
\frac{1}{\tau^{\beta}}\biggl[  \sum_{k=2}^{n} b_{n-k}^{(\beta)}(v^k-2v^{k-1}+v^{k-2})+ 2b_{n-1}^{(\beta)}(v^1-v^0) \biggl].
\end{eqnarray}  It holds that
\begin{eqnarray} \label{denotion-be}
&& _0^CD_t^{\beta}v(t_n)=\Delta_t^{\beta}v^n- 2  \frac{b_{n-1}^{(\beta)}}{\tau^{\beta-1}}  v^{\prime}(t_0)  +R_2[v(t_n)],\quad n\geq 1,
\end{eqnarray}
where
$$ |R_2[v(t_n)]|\leq \frac{9t_n^{2-\beta}}{\Gamma{(3-\beta)}} \max_{0\leq k\leq n-1}\int_{0}^{1}|v^{\prime\prime\prime}(t_{k}+\theta \tau)|d\theta \cdot \tau.$$
\end{lem}
\begin{proofed}
Note that
\begin{eqnarray}
&&_0^CD_t^{\beta}v(t_n)=\frac{1}{\Gamma{(2-\beta)} } \sum_{k=2}^n \int_{t_{k-1}}^{t_k} \frac{v^{\prime\prime}(s)}{(t_n-s)^{\beta-1}}ds +\frac{1}{\Gamma{(2-\beta)} } \int_{0}^{t_1} \frac{v^{\prime\prime}(s)}{(t_n-s)^{\beta-1}}ds\nonumber \\
&&   \qquad \qquad ~ = \frac{1}{\Gamma{(2-\beta)} } \sum_{k=2}^n \int_{t_{k-1}}^{t_k} \frac{\delta_t^2 v^{k-1}}{(t_n-s)^{\beta-1}}ds +\frac{1}{\Gamma{(2-\beta)} } \int_{t_{0}}^{t_1} \frac{\Delta_t^2v^{0}}{(t_n-s)^{\beta-1}}ds+R_2[v(t_n)] \nonumber \\
&&   \qquad \qquad~=\Delta_t^{\beta}v^n- 2  \frac{b_{n-1}^{(\beta)}}{\tau^{\beta-1}}  v^{\prime}(t_0)  +R_2[v(t_n)], \label{split-1}
\end{eqnarray}
where
$$\delta_t^2v^{k-1}=\frac{v(t_k)-2v(t_{k-1})+v(t_{k-2}) }{\tau^2},\quad \Delta_t^2v^{0}=2\frac{v(t_1)-v(t_{0})-\tau v^{\prime}(t_{0}) }{\tau^2},$$
and
\begin{eqnarray}
&&R_2[v(t_n)]= \frac{1}{\Gamma{(2-\beta)} } \sum_{k=2}^n \int_{t_{k-1}}^{t_k} \frac{v^{\prime\prime}(s)-\delta_t^2 v^{k-1}}{(t_n-s)^{\beta-1}}ds  +\frac{1}{\Gamma{(2-\beta)} }  \int_{t_{0}}^{t_1} \frac{v^{\prime\prime}(s)-\Delta_t^2 v^{0}}{(t_n-s)^{\beta-1}}ds. \label{truncation}
\end{eqnarray}
By  Taylor expansion with the integral remaining term,  we find for $s\in(t_{k-1},t_{k})$ that
\begin{eqnarray}
&& |v^{\prime\prime}(s)-\delta_t^2 v^{k-1}| = \biggl|\int^s_{t_{k-2}}\biggl(1-\frac{1}{2\tau^2}[(t_k-t)^2+(t_k-s)^2]\biggl) v^{\prime\prime\prime}(t)dt-\int_s^{t_{k}}\frac{1}{2\tau^2}(t_{k}-t)^2 v^{\prime\prime\prime}(t)dt\nonumber\\
&&\qquad\qquad\qquad\qquad +\int_{t_{k-2}}^{t_{k-1}}\frac{1}{\tau^2}(t_{k-1}-t)^2 v^{\prime\prime\prime}(t)dt \biggl|\nonumber\\
&& \qquad \qquad\qquad\quad \leq \frac{7}{2}\int^s_{t_{k-2}} |v^{\prime\prime\prime}(t)|dt+\frac{1}{2}\int_s^{t_{k}} |v^{\prime\prime\prime}(t)|dt+\int_{t_{k-2}}^{t_{k-1}} |v^{\prime\prime\prime}(t)|dt\nonumber\\
&& \qquad \qquad\qquad\quad \leq 9\max_{1\leq k\leq n}\int_{t_{k-1}}^{t_{k}} |v^{\prime\prime\prime}(t)|dt =9 \max_{0\leq k\leq n-1}\int_{0}^{1}|v^{\prime\prime\prime}(t_{k}+\theta \tau)|d\theta\cdot \tau . \label{truncation-1}
\end{eqnarray}
Similarly we can obtain
\begin{eqnarray}
 |v^{\prime\prime}(s)-\Delta_t^2 v^{0}|\leq \max_{0\leq k\leq n-1}\int_{0}^{1}|v^{\prime\prime\prime}(t_{k}+\theta \tau)|d\theta\cdot \tau. \label{truncation-2}
\end{eqnarray}
Combining  \eqref{truncation-1}-\eqref{truncation-2} with \eqref{truncation}, we get
\begin{eqnarray}
&& |R_2[v(t_n)]| \leq  \max_{0\leq k\leq n-1}\int_{0}^{1}|v^{\prime\prime\prime}(t_{k}+\theta \tau)|d\theta \cdot \frac{9}{\Gamma{(2-\beta)} } \sum_{k=1}^n \int_{t_{k-1}}^{t_k} \frac{1}{(t_n-s)^{\beta-1}}ds \cdot \tau
\nonumber \\
&&\qquad\qquad~
 = \frac{9t_n^{2-\beta}}{\Gamma{(3-\beta)}} \max_{0\leq k\leq n-1}\int_{0}^{1}|v^{\prime\prime\prime}(t_{k}+\theta \tau)|d\theta \cdot \tau.
\end{eqnarray}
This completes the proof.
\end{proofed}

\begin{rem}
When $\beta=2, $ \eqref{denotion-be} becomes the standard   first-order backward difference  approximation for second-order derivative, i.e.
$$v^{\prime\prime}(t_{n})=\frac{v(t_n)-2v(t_{n-1})+v(t_{n-2})}{\tau^2}+O(\tau),\quad n\geq2, $$
and
$$v^{\prime\prime}(t_{1})=\frac{2[v(t_1)-v(t_{0})-\tau v^{\prime}(t_{0})]}{\tau^2}+O(\tau). $$
\end{rem}

Take an integer $M.$ Let $\Omega_h\equiv \{ x_i~|~0\leq i\leq M\}$ be a uniform mesh of the interval $[0,L],$ where $x_i=ih,~0\leq i\leq M$ with $h=L/M.$ Suppose  $v=\{v_{i} \}$ is a grid function on $\Omega_h,$  define
$$\delta_xv_{i-\frac{1}{2}}=\frac{1}{h}(v_{i}-v_{i-1} ),\quad \delta_x^2 v_{i}=\frac{1}{h}(\delta_xv_{i+\frac{1}{2}}-\delta_xv_{i-\frac{1}{2}}). $$

To deal with spatial discretization and  construct a compact finite difference scheme for solving the problem \eqref{problem-1}-\eqref{problem-3}, we still need the lemma below.
\begin{lem}(See\cite{sun2005})\label{space-lem2}
Denote $\theta(s)=(1-s)^3[5- 3(1 - s)^2].$  If $w(x)\in  C^6[x_{i-1},x_{i+1}],$  it holds that
\begin{eqnarray*}
&&\frac{1}{12}[w^{\prime\prime}(x_{i-1})+10w^{\prime\prime}(x_i)+w^{\prime\prime}(x_{i+1})]
=\frac{w(x_{i-1})-2w(x_i)+w(x_{i+1})}{h^2}\\
&&\quad \quad \quad +\frac{h^4}{360}\int_0^1[w^{(6)}(x_i-sh)+w^{(6)}(x_i+sh)]\theta (s)ds,\quad 1\leq i\leq M-1.
\end{eqnarray*}
\end{lem}
The above Lemma can be  proved   straightforwardly by  Taylor expansion with  remaining integral.

Define  the average  operator $\mathcal{A}_x$ as
\begin{eqnarray*}
\mathcal{A}_xv_{i}=\left\{
\begin{array}{ll}
\frac{1}{12}(v_{i-1}+10v_{i}+v_{i+1}),~~~1\leq i\leq M-1,\\
v_{i},~~~~~~~~~~~~~~~~i=0~ \mbox{or}~ M.\\
\end{array}
\right.
\end{eqnarray*}
It is clear  that
$$\mathcal{A}_xv_{i}=(I+\frac{h^2}{12}\delta_x^2)v_{i},\quad 1\leq i\leq M-1, $$
where $I$ denotes the identity operator.
\subsection{Derivation of the difference scheme}
 We are now in a position to  derive the compact  difference scheme.
Define  grid functions below  $$U_i^n=u(x_i,t_n), \quad (U_t)_i^0= u_t(x_i,t_0)\quad F_i^n=f(x_i,t_n),\quad 0\leq i\leq M,\quad 0\leq n\leq N. $$
Considering  Equation \eqref{problem-1} at grid points  $ (x_i,t_n),$ we have
\begin{eqnarray}\label{eq-5}
&&K_1\  _0^{C}D_t^{\alpha}u(x_i,t_n)+K_2\  _0^{C}D_t^{\beta}u(x_i,t_n)=  \partial_x^2 u(x_i,t_n) +f(x_i,t_n),\quad 0\leq i\leq M,\quad 0\leq n\leq N.
\end{eqnarray}
 By Lemmas \ref{dis-lem-sun-1}, \ref{dis-lemma-2}, we have
\begin{eqnarray}\label{eq-6}
&&K_1 \delta_t^{\alpha}U_i^n  +K_2 \Delta_t^{\beta} U_i^n= \partial_x^2u (x_i,t_n)+F_i^n+\frac{2K_2b_{n-1}^{(\beta)}}{\tau^{\beta-1}} u_t(x_i,0)+K_1R_t^1[u(x_i,t_n)]+K_2R_t^2[u(x_i,t_n)], \nonumber\\
&&\quad 0\leq i\leq M,\quad 1\leq n\leq N.
\end{eqnarray}
Here $R_1$ and $R_2$ are similarly  defined as that in Lemmas \ref{dis-lem-sun-1} and \ref{dis-lemma-2} respectively.
For the spatial  discretization, acting the average operator $\mathcal{A}_x$, it follows from  Lemma \ref{space-lem2}  that
\begin{eqnarray}\label{eq-7}
&&K_1 \mathcal{A}_x\delta_t^{\alpha}U_i^n  +K_2 \mathcal{A}_x\Delta_t^{\beta} U_i^n=  \delta_x^2 U_i^n+\mathcal{A}_xF_i^n+   \frac{2K_2b_{n-1}^{(\beta)}}{\tau^{\beta-1}} \mathcal{A}_x (U_t)_i^0+R_i^n, \quad 1\leq i\leq M-1,\quad 1\leq n\leq N.\nonumber \\
\end{eqnarray}
where
$$R_i^n=K_1\mathcal{A}_xR_t^1[u(x_i,t_n)]+K_2\mathcal{A}_xR_t^2[u(x_i,t_n)]+R_x^3[u(x_i,t_n)].$$
 and $$R_x^3[u(x_i,t_n)]=\frac{h^4}{360}\int_0^1[\partial_x^6u(x_i-sh,t_n)+\partial_x^6u(x_i+sh,t_n)]\theta (s)ds.$$

It follows from Lemmas \ref{dis-lem-sun-1},  \ref{dis-lemma-2} and \ref{space-lem2}  that  there exists a constant  $C_u,$ which depends on the regularity of the solution $u(x,t)$ and the the parameters $\alpha$ and $\beta$ but is independent of the step size $h$ and $\tau,$
such that
\begin{eqnarray}
&&|R_i^n|\leq C_u(\tau+h^4),\quad 1\leq i\leq M-1,~ 1\leq n\leq N.\label{truncation-error-1}
\end{eqnarray}
Noticing the initial-boundary conditions,
\begin{eqnarray}
&&U_i^0=\phi_0(x_i) , \quad  (U_t)_i^0=\phi_1(x_i),\quad  0\leq i\leq M , \label{eq-8}   \\
&& U_0^n=\varphi_0(x_0),\quad U_M^n=\varphi_1(x_M),\quad 1\leq n\leq N ,     \label{eq-9}
\end{eqnarray}
 omitting the small terms $R_i^n $ and
denoting by $u_i^n$ the numerical approximation of $U_i^n,$  we get the compact finite difference scheme,
\begin{eqnarray}
&&K_1 \mathcal{A}_x\delta_t^{\alpha}u_i^n  +K_2 \mathcal{A}_x\Delta_t^{\beta} u_i^n=  \delta_x^2 u_i^n+\mathcal{A}_xF_i^n+   \frac{2K_2b_{n-1}^{(\beta)}}{\tau^{\beta-1}} \mathcal{A}_x \phi_i^1,
\quad
1\leq i\leq M-1,~1\leq n\leq  N,\quad \label{b9}\\
&&u_i^0=\phi_0(x_i) , \quad  0\leq i\leq M , \label{b10}   \\
&&  u_0^n=\varphi_0(t_n),\quad u_M^n=\varphi_1(t_n),\quad 1\leq n\leq N.  \label{b11}
\end{eqnarray}
\begin{rem}
When $\alpha=1,$ $\beta=2,$ we get the following three time-level backward difference scheme
\begin{eqnarray}
&&K_1 \mathcal{A}_x \biggl(\frac{u_i^n-u_i^{n-1}}{\tau}\biggl)  +K_2 \mathcal{A}_x \biggl( \frac{u_i^n-2u_i^{n-1}+u_i^{n-2}}{\tau^2} \biggl)=  \delta_x^2 u_i^n+\mathcal{A}_xF_i^n,
 \nonumber\\
&&
\quad
1\leq i\leq M-1,~2\leq n\leq  N,\quad \label{b99}\\
&&K_1 \mathcal{A}_x \biggl(\frac{u_i^1-u_i^{0}}{\tau} \biggl)  + 2K_2  \mathcal{A}_x \biggl(\frac{u_i^1-u_i^{0}-\tau \phi_i}{\tau^2} \biggl) =  \delta_x^2 u_i^n+\mathcal{A}_xF_i^n,
\quad
1\leq i\leq M-1,\quad \label{b112}\\
&&u_i^0=\phi_0(x_i) , \quad  0\leq i\leq M , \label{b100}   \\
&&  u_0^n=\varphi_0(t_n),\quad u_M^n=\varphi_1(t_n),\quad 1\leq n\leq N.  \label{b111}
\end{eqnarray}
\end{rem}

\section{Error analysis of the scheme \eqref{b9}-\eqref{b11}}
 In this section, we shall give  the stability and convergence analysis for the scheme \eqref{b9}-\eqref{b11}.
Let
$$V_h=\{v~|~ v ~\mbox{is a grid function on}~ \Omega_h ~\mbox{and } ~v_0=v_M=0\}.$$
For any $u, v\in V_h,$ we define the discrete inner products
$$(u,v)=h\sum_{i=1}^{M-1}u_iv_i,\qquad   \langle  \delta_xu,\delta_xv \rangle=h\sum_{i=1}^{M}\delta_xu_{i-1/2}\cdot \delta_xv_{i-1/2},   $$
and induced  norms
$$ ~ \|u\|=\sqrt{(u,u)},\qquad |u|_1=\sqrt{ \langle  \delta_xu,\delta_xu \rangle}.$$
Denote   maximum norm  by
$$\|u\|_{\infty}=\max_{0\leq i\leq M }|u_i|. $$
It is easy to check that
\begin{equation}\label{df-b2}
  (\delta_x^2u,v)=-\langle \delta_xu,\delta_xv \rangle.
\end{equation}
In addition, one has from \cite{sun2005} that
\begin{equation}
\|u\|_{\infty}\leq \frac{\sqrt{L}}{2} |u|_1. \label{dis-dif-b22}
\end{equation}

Some additional  lemmas are still required in order to prove the stability and  convergence of the proposed  scheme \eqref{b9}-\eqref{b11}.
\begin{lem}(See \cite{HaoLS2015}) \label{space-dis-lem-1}
For any $u,v\in V_h,$ there exists a positive definite operator  denoted by $ \mathcal{Q}_x$ such that
$$(\mathcal{A}_xu,v)=(\mathcal{Q}_xu,\mathcal{Q}_xv).$$
\end{lem}

\begin{lem}(See \cite{HaoLS2015})\label{lem5}
For any $v\in V_h,$ it holds that
$$\frac{2}{3}\|v\|^2 \leq (\mathcal{A}_xv,v)\leq \|v\|^2.$$
\end{lem}

Then we can define the equivalent weighted  norm  as
\begin{equation}\label{c18}
 \|v\|_A= \sqrt{(\mathcal{A}_xv,v)},
\end{equation}
and it follows that
\begin{equation}\label{eq:b3}
 \frac{2}{3}\|v\|^2 \leq \|v\|_A^2\leq \|v\|^2.
\end{equation}

\begin{lem}(See\cite{Marcos1990})\label{Lopez90}
Let $\{c_0, c_1,\ldots, c_n, \ldots \}$ be a sequence of real numbers with the properties below
$$c_n \geq 0,\quad c_n-c_{n-1} \leq 0, \quad c_{n+1}-2c_n +c_{n-1}\geq 0. $$
Then for any positive integer $m$ and  each vector $(v_1,v_2,\ldots, v_m)$ with $m$ real entries, it holds that
$$\sum_{n=1}^m\biggl(\sum_{p=0}^{n-1}c_pv_{n-p}\biggl)v_n \geq 0. $$
\end{lem}
Denote $\delta_tv^n=\frac{1}{\tau}(v^n-v^{n-1}).$
Notice that  $\delta_t^{\alpha} v^n$ can be reformulated as
$$\delta_t^{\alpha} v^n=\tau^{1-\alpha} \sum_{k=0 }^{n-1} a_{k}^{(\alpha)}\delta_tv^{n-k}.  $$  It is not difficult to verify that the   coefficients $\{a_k^{(\alpha)} \}$ defined by \eqref{coefficients}   satisfy
\begin{eqnarray}
&&1=a_0^{(\alpha)} > a_1^{(\alpha)} > a_2^{(\alpha)}>\cdots  > a_k^{(\alpha)} >\cdots \rightarrow 0,\label{HL-2}\\
&& \frac{(k+1)^{-\alpha}}{\Gamma{(1-\alpha)}}<  a_k^{(\alpha)} < \frac{k^{-\alpha}}{\Gamma{(1-\alpha)}}\\
&& a^{(\alpha)}_{k+1}-2a^{(\alpha)}_k+a^{(\alpha)}_{k-1} \geq 0.
\end{eqnarray}

Thus, by Lemma  \ref{Lopez90},  we have
\begin{eqnarray}
\sum_{n=1}^m(\delta_t^{\alpha} v^n)\cdot  (\delta_t v^n)=\tau^{1-\alpha} \sum_{n=1}^m\biggl(\sum_{k=0 }^{n-1} a_{k}^{(\alpha)}\delta_tv^{n-k}\biggl) \cdot (\delta_t v^n) \geq 0.\label{Ineq:semipositive}
\end{eqnarray}

\begin{lem}\label{Lem-embedding-1}
For any $v=\{v^0,v^1,v^2,v^3,\ldots  \}$  and $\alpha \in (0,1),$ we have
\begin{eqnarray}\label{dis-positivety}
\tau \sum_{n=1}^{m} v^n\delta_t^{\alpha}  v^n  \geq \frac{1}{ 2} \tau^{1-\alpha} \sum_{n=1}^{m}a_{m-n}^{(\alpha)}(v^n)^2-\frac{t_m^{1-\alpha}}{2\Gamma(2-\alpha)}(v^0)^2.
\end{eqnarray}
\end{lem}
\begin{proofed}
Following the proof technique of \cite{Alikhanov2015}, we can prove that there holds  the discrete analogy of the energy inequality i.e.,
\begin{equation}\label{energy:ineq-dis}
v^n\delta_t^{\alpha}v^n \geq \frac{1}{2} \delta_t^{\alpha}(v^n)^2.
\end{equation}
Thus it suffices to  prove the following inequality
\begin{equation}\label{energy:ineq-dis-3}
 \frac{1}{2} \tau\sum_{n=1}^m\delta_t^{\alpha}(v^n)^2\geq \frac{1}{ 2} \tau^{1-\alpha} \sum_{n=1}^{m}a_{m-n}^{(\alpha)}(v^n)^2-\frac{t_m^{1-\alpha}}{2\Gamma(2-\alpha)}(v^0)^2.
\end{equation}
 On the one hand, one has
\begin{eqnarray*}
&&\frac{1}{2} \tau\sum_{n=1}^m\delta_t^{\alpha}(v^n)^2=\frac{1}{2} \tau^{2-\alpha}\sum_{n=1}^m \sum_{k=1}^{n} a_{n-k}^{(\alpha)}\delta_t(v^k)^2=\frac{1}{2} \tau^{2-\alpha}\sum_{k=1}^m\delta_t(v^k)^2 \sum_{n=k}^{m} a_{n-k}^{(\alpha)}\\
 &&\quad=\frac{1}{2} \tau^{1-\alpha}\sum_{k=1}^m[(v^k)^2 -(v^{k-1})^2]  \frac{(m-k+1)^{1-\alpha}}{\Gamma(2-\alpha)}.
\end{eqnarray*}
On the other hand
\begin{eqnarray*}
  &&\quad\frac{1}{2} \tau^{1-\alpha}\sum_{k=1}^m[(v^k)^2 -(v^{k-1})^2]  \frac{(m-k+1)^{1-\alpha}}{\Gamma(2-\alpha)}\\
    &&=\frac{1}{2\Gamma(2-\alpha)} \tau^{1-\alpha}\biggl(\sum_{k=1}^{m}(v^k)^2[ (m-k+1)^{1-\alpha}-(m-k)^{1-\alpha}]  -(v^0)^2m^{1-\alpha}\biggl)\\
    &&=\frac{1}{2} \tau^{1-\alpha}\sum_{k=1}^{m}a_{n-k}^{(\alpha)}(v^k)^2 -\frac{t_m^{1-\alpha}}{2\Gamma(2-\alpha)}(v^0)^2,
\end{eqnarray*}
This concludes the proof of Lemma \ref{Lem-embedding-1}.
\end{proofed}

It should be noted that Sun and Wu also  gave the proof of the equality \eqref{dis-positivety} in \cite{Sun2006}. But here we adopted a different technique  in \cite{Sun2006}.

We now turn to investigate  the convergence of the  scheme \eqref{b9}-\eqref{b11}.

\begin{thm}[Convergence]\label{Convergence}
Suppose  $u(x,t) $ solves the problem \eqref{problem-1}-\eqref{problem-3}, $ u(x,t)\in C_{x,t}^{5,2}([0,L]\times [0,T])$ and   $ \partial_x^6u(x,\cdot)\in L^1[0,L]$ $\partial_t^3u(\cdot,t)\in  L^1[0,T].$  Let $\{u_i^n\}$ be the solution of difference scheme \eqref{b9}-\eqref{b11}.
Then for $n\tau \leq T,$ there exists a constant $C$ such that
$$\max_{0\leq i\leq M}|u(x_i,t_n)-u_i^n|\leq C(
\tau+h^4). $$
\end{thm}
\begin{proofed}
Let $$e_i^n=U_i^n-u_i^n,\quad 0\leq i\leq M,~ 0\leq n\leq N.$$
Subtracting \eqref{b9}-\eqref{b11} from  \eqref{eq-7},\eqref{eq-8}-\eqref{eq-9}, we get the error equations,
\begin{eqnarray}
&&K_1 \mathcal{A}_x\delta_t^{\alpha}e_i^n  +K_2 \mathcal{A}_x\Delta_t^{\beta} e_i^n=  \delta_x^2 e_i^n+R_i^n, \quad 1\leq i\leq M-1,~1\leq n\leq  N,\quad \label{dis-dif-b19}\\  
&&e_0^n=0,\quad e_M^n=0,\quad 1\leq n\leq N, \label{dis-dif-b20}   \\
&&e_i^0=0 , \quad  0\leq i\leq M   .     \label{dis-dif-b21}
\end{eqnarray}
Recall that
\begin{eqnarray*}
&&\Delta_t^{\beta}e^n=
\frac{1}{\tau^{\beta}}\biggl[  \sum_{k=2}^{n-2} b_{n-k}^{(\beta)}(e^k-2e^{k-1}+e^{k-2})+ 2b_{n-1}^{(\beta)}(e^1-e^0) \biggl].
\end{eqnarray*}
Denote $\delta_t e^n=\frac{1}{\tau}(e^n-e^{n-1}).$ Hereafter we stipulate $\delta_t e^0=0.$ Since $e^0=0,$ we can recast $\Delta_t^{\beta}e^n$ as
\begin{eqnarray}\label{HL-1}
&&\Delta_t^{\beta}e^n =\frac{1}{\tau^{\beta-1}}[b_0^{(\beta)}\delta_te^n-\sum_{k=1}^{n-1}(b_{n-k-1}^{(\beta)}-b_{n-k}^{(\beta)})\delta_te^k-b_{n-1}^{(\beta)}\delta_t e^0]+  \frac{ b_{n-1}^{(\beta)}}{\tau^{\beta}}e^1\nonumber\\
&&\qquad ~~=\delta_t^{\beta-1}(\delta_t e^n)+  \frac{ b_{n-1}^{(\beta)}}{\tau^{\beta}}e^1,
\end{eqnarray}
where we used the notation  \eqref{eq:def-al}.
Taking the inner product of \eqref{dis-dif-b19} with $\tau \delta_te^{n}$ and summing up for $ n=1$ to $m$ yield
\begin{eqnarray}\label{dis-dif-b12}
&& K_1 \tau \sum_{n=1}^{m}(\mathcal{A}_x\delta_t^{\alpha}e^n, \delta_te^{n} )  +K_2\tau \sum_{n=1}^{m} (\mathcal{A}_x\Delta_t^{\beta} e^n,\delta_te^{n})- \tau \sum_{n=1}^{m} (\delta_x^2 e_i^n,\delta_te^{n})=\tau \sum_{n=1}^{m}(R^n,\delta_te^{n}).
\end{eqnarray}
For the first  term on the left hand side of \eqref{dis-dif-b12}, using Lemma \ref{space-dis-lem-1} and  noting the commutation of the operators in different direction, and  the inequality \eqref{Ineq:semipositive}, we have
\begin{eqnarray} \label{dis-dif-b13-1}
&&   K_1 \tau \sum_{n=1}^{m}(\mathcal{A}_x\delta_t^{\alpha}e^n, \delta_te^{n} )=  K_1 \tau \sum_{n=1}^{m}(\mathcal{Q}_x\delta_t^{\alpha}e^n,\mathcal{Q}_x\delta_te^{n} )=K_1 \tau\sum_{i=1}^{M} \sum_{n=1}^{m}(\delta_t^{\alpha}\mathcal{Q}_xe_i^n)( \delta_t \mathcal{Q}_xe_i^{n}) \geq 0  .
\end{eqnarray}
For the second term on the  left hand side of \eqref{dis-dif-b12}, by Lemmas \ref{space-dis-lem-1}, \ref{Lem-embedding-1} and  Equality  \eqref{HL-1},  we have
\begin{eqnarray} \label{dis-dif-b13-2}
&&K_2\tau \sum_{n=1}^{m} (\mathcal{A}_x\Delta_t^{\beta} e^n,\delta_te^{n})=K_2\tau \sum_{n=1}^{m} (\Delta_t^{\beta}\mathcal{Q}_x e^n,\mathcal{Q}_x\delta_te^{n}) \nonumber\\
&&=K_2\tau \sum_{n=1}^{m} (\delta_t^{\beta-1} \delta_t\mathcal{Q}_x e^n,\delta_t\mathcal{Q}_xe^{n})+K_2\tau \sum_{n=1}^{m}\frac{ b_{n-1}^{(\beta)}}{\tau^{\beta-1}}( \frac{\mathcal{Q}_x e^1}{\tau},\delta_t\mathcal{Q}_xe^{n}) \nonumber\\
&& \geq  \frac{K_2t_m^{1-\beta}}{ 2\Gamma(2-\beta)} \tau \sum_{n=1}^{m}\|\delta_t\mathcal{Q}_xe^{n}\|^2
-K_2\tau \sum_{n=1}^{m}\frac{ b_{n-1}^{(\beta)}}{\tau^{\beta-1}}  \biggl(\frac{\epsilon_1}{4}\| \frac{\mathcal{Q}_x e^1}{\tau}\|^2+ \frac{1}{\epsilon_1}\|\delta_t\mathcal{Q}_xe^{n}\|^2 \biggl)  \nonumber\\
&& \geq  \frac{K_2t_m^{1-\beta}}{ 2\Gamma(2-\beta)} \tau \sum_{n=1}^{m}\|\delta_t\mathcal{Q}_xe^{n}\|^2
-  \frac{\epsilon_1K_2 t_m^{2-\beta}}{4\Gamma(3-\beta)}\| \frac{\mathcal{Q}_x e^1}{\tau}\|^2-\frac{K_2}{\epsilon_1} \max_{0\leq n\leq m-1}\{\frac{ b_{n}^{(\beta)}}{\tau^{\beta-1}} \} \tau\sum_{n=1}^{m} \|\delta_t\mathcal{Q}_xe^{n}\|^2  \nonumber\\
&& = \biggl( \frac{K_2t_m^{1-\beta}}{ 2\Gamma(2-\beta)}-\frac{K_2}{\epsilon_1} \frac{ b_{0}^{(\beta)}}{\tau^{\beta-1}}\biggl) \tau \sum_{n=1}^{m}\|\delta_te^{n}\|_{A}^2
- \frac{\epsilon_1K_2  t_m^{2-\beta}}{4\Gamma(3-\beta)}\| \frac{ e^1}{\tau}\|_A^2.
\end{eqnarray}
where we have  used the $\epsilon$-type inequality and the  commutativity of the operators, i.e. $\mathcal{Q}_x\delta_t=\delta_t\mathcal{Q}_x.$

For the third   term on the left  hand side of \eqref{dis-dif-b12}, observing the initial value \eqref{dis-dif-b21}, we obtain
\begin{eqnarray}\label{dis-dif-b14}
&&- \tau \sum_{n=1}^{m} ( \delta_x^2 e^{n} , \delta_t e^{n} )=\tau \sum_{n=1}^{m} \langle  \delta_x e^{n} , \delta_t\delta_x  e^{n} \rangle = \sum_{n=1}^{m} \langle  \delta_x e^{n} , \delta_x  e^{n}-\delta_x  e^{n-1} \rangle  \nonumber \\
&&\geq \frac{1}{2} \sum_{n=1}^{m}(|  e^n|_1^2-| e^{n-1}|_1^2 )= \frac{1}{2} |  e^m|_1^2,
\end{eqnarray}
where  we have used the inequality $a(a-b)\geq \frac{1}{2}(a^2-b^2). $

For the   term on the right  hand side, we have
\begin{eqnarray}\label{dis-dif-b15}
&& \quad \tau \sum_{n=1}^{m} (  R^{n} , \delta_t e^{n} )  \leq  \tau \sum_{n=1}^{m}  \biggl (  \epsilon_2  \| \delta_t e^{n}\|^2 + \frac{1}{4\epsilon_2 } \|R^{n}\|^2 \biggl ).
 \end{eqnarray}
Substituting \eqref{dis-dif-b13-1}-\eqref{dis-dif-b15} into \eqref{dis-dif-b12}  gives
\begin{eqnarray*}
&&\quad  \biggl( \frac{K_2t_m^{1-\beta}}{ 2\Gamma(2-\beta)}-\frac{K_2}{\epsilon_1} \frac{ b_{0}^{(\beta)}}{\tau^{\beta-1}}\biggl) \tau \sum_{n=1}^{m}\|\delta_te^{n}\|_{A}^2
+\frac{1}{2}|e^m|_1^2 \nonumber \\
&&\leq\epsilon_2 \tau \sum_{n=1}^{m}\|\delta_te^{n}\|+    \frac{\epsilon_1 K_2 t_m^{2-\beta}}{4\Gamma(3-\beta)}\| \frac{ e^1}{\tau}\|_A^2+   \frac{1}{4\epsilon_2} \tau\sum_{n=1}^{m}    \|R^{n}\|^2  ,\quad 1\leq m \leq N.
 \end{eqnarray*}
 Noting the  norms  equivalence \eqref{eq:b3},  we have
 \begin{eqnarray*}
&&\quad  \biggl( \frac{K_2t_m^{1-\beta}}{ 3\Gamma(2-\beta)}-\frac{2K_2}{3\epsilon_1} \frac{ b_{0}^{(\beta)}}{\tau^{\beta-1}}\biggl) \tau \sum_{n=1}^{m}\|\delta_te^{n}\|^2+\frac{1}{2}|e^m|^2 \nonumber \\
&&\leq  \epsilon_2 \tau \sum_{n=1}^{m}\|\delta_tv^{n}\|+     K_2  \frac{\epsilon_1 t_m^{2-\beta}}{4\Gamma(3-\beta)}\| \frac{ e^1}{\tau}\|_A^2+   \frac{1}{4\epsilon_2} \tau\sum_{n=1}^{m}    \|R^{n}\|^2  ,\quad 1\leq m \leq N.
 \end{eqnarray*}
Recall that $b_0^{(\beta)}=a_0^{(\beta-1)}=\frac{1}{\Gamma(3-\beta)}.$
Taking $\epsilon_1=\frac{4\tau^{(1-\beta)}t_m^{\beta-1}}{2-\beta}$, $\epsilon_2=\frac{K_2t_m^{1-\beta}}{ 12\Gamma(2-\beta)}$  leads to
 \begin{eqnarray}\label{HL-4}
&&\quad   \frac{K_2t_m^{1-\beta}}{ 12\Gamma(2-\beta)} \tau \sum_{n=1}^{m}\|\delta_te^{n}\|^2+\frac{1}{2}|e^m|^2 \leq       \frac{  K_2 t_m}{(2-\beta)\Gamma(3-\beta)\tau^{1+\beta}}\| e^1\|_A^2+   \frac{ 3\Gamma(2-\beta)t_m^{\beta-1}}{K_2} \tau\sum_{n=1}^{m}    \|R^{n}\|^2  ,\nonumber\\
&&\quad 1\leq m \leq N.
 \end{eqnarray}

 Taking the inner product of \eqref{dis-dif-b19} with $ \tau\delta_te^{1}$ gives
\begin{eqnarray*}
&& K_1 \tau(\mathcal{A}_x\delta_t^{\alpha}e^1, \delta_te^{1} )  +K_2\tau (\mathcal{A}_x\Delta_t^{\beta} e^1,\delta_te^{1})- \tau (\delta_x^2 e_i^1,\delta_te^{1})=\tau(R^1,\delta_te^{1}).
\end{eqnarray*}
That is
\begin{eqnarray*}
&& \frac{K_1 a_0^{(\alpha)}}{\tau^{\alpha}}(\mathcal{A}_xe^1 , e^1 )  + \frac{2K_2 b_0^{(\beta)}}{\tau^{\beta}}(\mathcal{A}_xe^1 , e^1)-  (\delta_x^2 e^1, e^1)=(R^1, e^1).
\end{eqnarray*}
Similar to the argument as above, we can obtain
\begin{eqnarray*}
&&  \frac{2K_2 b_0^{(\beta)}}{\tau^{\beta}}(\mathcal{A}_xe^1 , e^1)\leq (R^1, e^1)\leq \|R^1\|\cdot \|e^1\| \leq \|R^1\|\cdot  \frac{3}{2}\|e^1\|_A.
\end{eqnarray*}
Thus, we have
\begin{eqnarray}\label{HL-3}
&&  \|e^1\|_{A} \leq \frac{3\tau^{\beta}}{4K_2 b_0^{(\beta)}}\|R^1\|=\frac{3\Gamma(3-\beta)\tau^{\beta}}{4K_2 }\|R^1\|  .
\end{eqnarray}
Inserting \eqref{HL-3} into \eqref{HL-4} leads to
 \begin{eqnarray}\label{HL-5}
&&\quad   \frac{K_2t_m^{1-\beta}}{ 12\Gamma(2-\beta)} \tau \sum_{n=1}^{m}\|\delta_te^{n}\|^2+\frac{1}{2}|e^m|^2 \nonumber\\
&& \leq       \frac{  9\Gamma(2-\beta) t_m\tau^{\beta-1}}{16K_2}\| R^1\|^2+   \frac{ 3\Gamma(2-\beta)t_m^{\beta-1}}{K_2} \tau\sum_{n=1}^{m}    \|R^{n}\|^2  \nonumber\\
&& \leq \frac{3\Gamma(2-\beta)t_m}{K_2}(\frac{3}{4}\tau^{\beta-1}+t_m^{\beta-1} )LC_u^2(\tau+h^4)^2,\quad 1\leq m \leq N.
 \end{eqnarray}
It follows from \eqref{dis-dif-b22}  that
\begin{eqnarray}
&&\|e^m\|_{\infty}\leq \frac{\sqrt{L}}{2}  |e^m|_1 \leq \sqrt{\frac{3\Gamma(2-\beta)}{K_2}(\frac{3}{4}\tau^{\beta-1}+t_m^{\beta-1} )}LC_ut_m(\tau+h^4):= C(\tau+h^4) .\label{ineq:error-1}
 \end{eqnarray}
This completes the proof.
\end{proofed}

With the above  convergency  and consistency condition \eqref{truncation-error-1}, the stability of the difference scheme can be  obtained  straightforwardly.
 \begin{cor}[Stability]
The difference scheme \eqref{b9}-\eqref{b11} is unconditionally stable to the initial time value and    the right hand term.
\end{cor}


\section{Extension to the distributed order equation}
  The diffusion-wave equation with distributed order in time can be written as:
\begin{eqnarray}\label{HL-6}
\int_a^bw(\alpha)_0^{C}D_t^{\alpha}u(x,t)d\alpha=\partial_x^{2}u(x,t)+f(x,t), \quad 0<x<L, \quad t>0,
\end{eqnarray}
where $0<a<b<2,$ and the function $w(\alpha)$ acting as weight for the order of
differentiation is subject to that $w(\alpha)>0 $ and $\int_{a}^bw(\alpha)=Constant>0.$ If $a = 0$ and $ b = 1$ we obtain the sub-diffusion
equation with distributed order in time, and if $a = 1$ and $b = 2$, we have the diffusion-wave equation with distributed order in time. As an extension of multi-term equations, distributed order equations have gained considerable attention   recently. For the non-Markovian process which is non self-similar and exhibits a continuous distribution of time-scales, a continuous
distribution of fractional time derivatives is introduced in \cite{Mainardi2007}. For example, some complicated processes involving a mixture of
power laws often lead to the distributed-order fractional derivative in time \cite{Chechkin2002,Meerschaert2011}.   As a
precise tool to explain and describe some real physical phenomena, numerical work of the time fractional differential equations involved with the distributed order operator has also attract the attention of many scholars. For example,
Diethelm and Ford \cite{Diethelm2001} introduced a general framework for distributed-order ordinary
differential equations by using the quadrature formula, such as the trapezoidal formula, with
some suitable numerical solver for the resulting multi-term fractional equations, while a convergence
analysis of the method was discussed in \cite{Diethelm2009} recently.
Numerical schemes for partial integro-differential equations with distributed fractional order,
including the  sub-diffusion equation \eqref{HL-6}, have  appeared in the literature very recently, see \cite{GaoSun2015,Ye2015,Morgato2015}. In this work, we are concerned with the case $a=0, b=2,$ which, to the best of our knowledge, has not been investigated to date.

In this section, the implicit difference  scheme \eqref{b9}-\eqref{b11} is extended to approximate the distributed-order diffusion-wave equation \eqref{HL-6} with $a=0, b=2$. The initial boundary conditions \eqref{problem-2}-\eqref{problem-3}
 are considered for equation \eqref{HL-6}.   Given a positive integer $J,$ take $\sigma=\frac{1}{J}.$ We partition the interval $[0,2]$ as follows:
$$0=\beta_0<\beta_1<\ldots<\beta_J=1<\beta_{J+1}<\ldots<\beta_{2J}=2,$$
where $\beta_{j}=j\sigma$ for $j=0,1,2,\cdots,2J.$
We can use the mid-point quadrature rule for approximating the integral in \eqref{HL-6}. Let $\alpha_j=\frac{\beta_{j-1}+\beta_j}{2}.$  Then,
\begin{equation}\label{HL-7}
\int_0^2w(\alpha)_0^CD_t^{\alpha}u(x,t)d\alpha=\sigma\sum_{j=1}^{2J}K_j\  _0^CD_t^{\alpha_j}u(x,t)-
\frac{\sigma^2}{24} \Phi^{\prime\prime}(\zeta),
\end{equation}
where $\zeta\in (0,2),$ $\Phi(\alpha)=w(\alpha) _0^CD_t^{\alpha}u(x,t)$ and $K_j=w(\alpha_j).$
By the formula \eqref{HL-7}, the distributed order fractional diffusion-wave equation are reduced to the following multi-term TFMDWEs as follows:
\begin{eqnarray}\label{HL-8}
\sigma\sum_{j=1}^{2J}K_j \ _0^CD_t^{\alpha_j}u(x,t)=\partial_x^{2}u(x,t)+f(x,t)+\frac{\sigma^2}{24} \Phi^{\prime\prime}(\zeta),\quad t>0, \quad 0<x<L.
\end{eqnarray}
Considering  the equation \eqref{HL-8} at grid points  $ (x_i,t_n),$ we have
\begin{eqnarray}\label{HL-9}
&&\sigma\sum_{j=1}^{J}K_j\  _0^{C}D_t^{\alpha_j}u(x_i,t_n)+\sigma\sum_{j=J+1}^{2J}K_j\  _0^{C}D_t^{\alpha_j}u(x_i,t_n)=  \partial_x^2 u(x_i,t_n) +f(x_i,t_n)+\frac{\sigma^2}{24} \Phi^{\prime\prime}(\zeta),\quad 0\leq i\leq M,\nonumber\\
&&\quad 1\leq n\leq N.
\end{eqnarray}
Following the similar  derivation as that in  Subsection 3.2, we can obtain
\begin{eqnarray}\label{HL-10}
&&\sigma\sum_{j=1}^{J}K_j \mathcal{A}_x\delta_t^{\alpha_j}U_i^n  +\sigma\sum_{j=J+1}^{2J}K_j \mathcal{A}_x\Delta_t^{\alpha_j} U_i^n=  \delta_x^2 U_i^n+\mathcal{A}_xF_i^n+ \sigma\sum_{j=J+1}^{2J}K_j  \frac{2b_{n-1}^{(\alpha_j)}}{\tau^{\alpha_j-1}} \mathcal{A}_x (U_t)_i^0+\tilde{R}_i^n, \nonumber\\
&&\quad 1\leq i\leq M-1,\quad 1\leq n\leq N,
\end{eqnarray}
where
\begin{equation}\label{HL-11}
|\tilde{R}_i^n|\leq \tilde{C}_u(\tau+\sigma^2+h^4).
\end{equation}
Thus, we get the following implicit  compact finite difference scheme
\begin{eqnarray}
&&\sigma\sum_{j=1}^{J}K_j \mathcal{A}_x\delta_t^{\alpha_j}u_i^n  +\sigma\sum_{j=J+1}^{2J}K_j \mathcal{A}_x\Delta_t^{\alpha_j} u_i^n=  \delta_x^2 u_i^n+\mathcal{A}_xF_i^n+   \sigma\sum_{j=J+1}^{2J}K_j  \frac{2b_{n-1}^{(\alpha_j)}}{\tau^{\alpha_j-1}} \mathcal{A}_x \phi_i^1,
 \nonumber\\
&&
\quad
1\leq i\leq M-1,~1\leq n\leq  N,\quad \label{HL-12}\\
&&u_i^0=\phi_0(x_i) , \quad  0\leq i\leq M , \label{HL-13}   \\
 &&u_0^n=\varphi_0(t_n),\quad u_M^n=\varphi_1(t_n),\quad 1\leq n\leq N .  \label{HL-14}
\end{eqnarray}

Following similar arguments of Section 4, we have the following result.
\begin{thm}[convergence and stability]\label{Convergence and Stability}
The scheme \eqref{HL-12}-\eqref{HL-14} is unconditionally stable. Furthermore,  let $w(\alpha)\in C[0,2]$ and $\Phi(\alpha)\in C^2[0,2].$ suppose  $u(x,t) $ solves the problem \eqref{HL-6} with initial-boundary value \eqref{problem-2}-\eqref{problem-3}, $ u(x,t)\in C_{x,t}^{5,2}([0,L]\times [0,T])$ and   $ \partial_x^6u(x,\cdot)\in L^1[0,L],$ $\partial_t^3u(\cdot,t)\in  L^1[0,T].$  And $\{u_i^n\}$ be the solution of difference scheme \eqref{HL-12}-\eqref{HL-14}.
Then for $n\tau \leq T,$ there exists a constant $C$ such that
$$\max_{0\leq i\leq M}|u(x_i,t_n)-u_i^n|\leq C(
\tau+\sigma^2+h^4). $$
\end{thm}
It should be noted that Gauss quadrature can also be used to approximate the integration with high accuracy when the  integrand $\Phi(\alpha)$ is smooth enough. The subsequent  procedure is as identical as the mid-point case, here we skip it.
\section{Extension to the two-dimensional case}
In this section, we will consider the generalization of the our proposed method to the following two-dimensional equation
\begin{eqnarray}\label{multi-dw-2D}
&&K_1\  _0^{C}D_t^{\alpha}u(x,y,t)+K_2\  _0^{C}D_t^{\beta}u(x,y,t)= \partial_x^{2}u(x,y,t)+\partial_y^{2}u(x,y,t) +f(x,y,t),\quad (x,y)\in \Omega,~t\in(0, T],\nonumber\\
\end{eqnarray}
which is subject to the initial-boundary conditions
\begin{eqnarray}
&&u(x,y,0)=\phi_0(x,y),\quad u_t(x,y,0)=\phi_1(x,y),\quad (x,y)\in \Omega,\label{initial-2D} \\
 &&u(x,y,t)=0,\quad (x,y)\in \partial{\Omega}, \quad  t\in(0, T].\label{boundary-2D}
 \end{eqnarray}
 where $\Omega=(0,L_1)\times(0,L_2).$
For the  spatial approximation, take  two integers $M_1,M_2 $ and let
$h_1=(b-a)/M_1,$ $h_2=(d-c)/M_2,$
 $x_i=ih_1,$ $0\leq i\leq M_1,$ $y_j=j h_2, $ $0\leq j\leq M_2.$
Let $\bar{\Omega}_h=\{(x_i,y_j)|0\leq i\leq M_1,~ 0\leq j\leq M_2\},$ and $\Omega_h=\bar{\Omega}_h\cap\Omega, $ and $\partial \Omega_h=\bar{\Omega}_h\cap\partial\Omega. $
For any grid function $v=\{v_{i,j}|0\leq i\leq M_1, 0\leq j\leq M_2\},$ denote
$$\delta_xv_{i-\frac{1}{2},j}=\frac{1}{h_1}(v_{i,j}-v_{i-1,j} ),\quad \delta_x^2 v_{i,j}=\frac{1}{h_1}(\delta_xv_{i+\frac{1}{2},j}-\delta_xv_{i-\frac{1}{2},j}). $$
Similar notations $\delta_yv_{i,j-\frac{1}{2}}, $ $\delta_y^2v_{i,j} $ can be defined.
The spatial average  operators are defined as
\begin{eqnarray*}
\mathcal{A}_xv_{i,j}=\left\{
\begin{array}{ll}
\frac{1}{12}(v_{i-1,j}+10v_{i,j}+v_{i+1,j}),~~~1\leq i\leq M_1-1,~~0\leq j\leq M_2,\\
v_{i,j},~~~~~~~~~~~~~~~~i=0~ \mbox{or}~ M_1,~~0\leq j\leq M_2,\\
\end{array}
\right.\\
\mathcal{A}_yv_{i,j}=\left\{
\begin{array}{ll}
\frac{1}{12}(v_{i,j-1}+10v_{i,j}+v_{i,j+1}),~~~1\leq j\leq  M_2-1,~~0\leq i\leq M_1,\\
v_{i,j},~~~~~~~~~~~~~~~~j=0 ~\mbox{or} ~M_2,~~0\leq i\leq M_1.\\
\end{array}
\right.
\end{eqnarray*}
Akin to the construction as the one-dimensional case,  we present the  implicit   compact difference scheme  for  the problem  \eqref{multi-dw-2D}-\eqref{boundary-2D} as follows
\begin{eqnarray}
&&K_1\mathcal{A}_y \mathcal{A}_x\delta_t^{\alpha}u_{i,j}^n  +K_2\mathcal{A}_y  \mathcal{A}_x\Delta_t^{\beta} u_{i,j}^n= \mathcal{A}_y  \delta_x^2 u_{i,j}^n+\mathcal{A}_x  \delta_y^2 u_{i,j}^n+\mathcal{A}_y\mathcal{A}_xF_{i,j}^n+   \frac{2K_2b_{n-1}^{(\beta)}}{\tau^{\beta-1}} \mathcal{A}_y \mathcal{A}_x \phi_{i,j}^1,
 \nonumber\\
&&
\quad
(x_i,y_j)\in \Omega_h,~1\leq n\leq  N,\quad \label{b9-2D}\\
&&u_{i,j}^0=\phi_0(x_i,y_j) , \quad (x_i,y_j)\in \Omega_h, \label{b10-2D}   \\
&&  u_{i,j}^n=0,\quad (x_i,y_j)\in \Omega_h,~1\leq n\leq N .  \label{b11-2D}
\end{eqnarray}
  The convergence in  $H_1$ norm and stability can be derived by the similar analysis as one-dimensional case.  Here, we will not dwell on this anymore.

\section{Implementation of scheme \eqref{b9}-\eqref{b11}}
\setcounter{equation}{0}
In this section, we propose a fast solver based on the partial diagonalization
 technique together with dived and conquer strategy. To fix the idea, we consider the homogeneous boundary condition for the simplicity.  Several lemmas are  presented first.
\begin{lem}\label{lem-trimax-1} \cite{ZhangD2014}
A general tri-diagonal Toeplitz matrix of order $n-1$ is given as
 \begin{eqnarray*}
 \mathbf{T}=\left(\begin{array}{ccccccc}
b&c&0&\cdots&0&0\\
a&b&c&\cdots&0&0\\
0&a&b&\cdots&0&0\\
\vdots&\vdots&\vdots&\ddots&\vdots&\vdots\\
0&0&0&\cdots&b&c\\
0&0&0&\cdots&a&b\\
\end{array}\right).
\end{eqnarray*}
the eigenvalues and eigenvectors of the tri-diagonal Toeplitz matrix $\mathbf{T}$ are given by
$$\lambda_i=b+2a\sqrt{\frac{c}{a}}\cos(\frac{\pi i}{n}),\quad i=1,2,\ldots,n-1,$$ and
 \begin{eqnarray*}
 \xi_i=\left(\begin{array}{c}
(\frac{c}{a})^{\frac{1}{2}}\sin(\frac{\pi}{n}i)\\
(\frac{c}{a})^{\frac{2}{2}}\sin(\frac{2\pi}{n}i)\\
(\frac{c}{a})^{\frac{3}{2}}\sin(\frac{3\pi}{n}i)\\
\vdots\\
(\frac{c}{a})^{\frac{n-1}{2}}\sin(\frac{(n-1)\pi}{n}i)\\
\end{array}\right),\quad i=1,2,\ldots n-1,
\end{eqnarray*}
i.e., $\mathbf{T}\xi_i=\lambda_i\xi_i,\quad i=1,2,\ldots n-1.$ The matrix $\mathbf{T}$ is diagonalizable and $\mathbf{P} = (\xi_1, \xi_2, \xi_3, \ldots, \xi_{n-1})$ diagonalizes $\mathbf{T}$, i.e., $ \mathbf{P}^{-1}\mathbf{T}\mathbf{P}=\mathbf{\Lambda},$
 where $\mathbf{\Lambda}=diag(\lambda_1,\lambda_2,\ldots,\lambda_{n-1}).$
 Moreover  $\mathbf{P}$ can be orthogonal, that is,$ \sqrt{\frac{2}{n}}\mathbf{P}$ is the orthogonal matrix.
\end{lem}

Now, we introduce the following  matrices:
 \begin{eqnarray*}
 \mathbf{M}_x=\left(\begin{array}{cccccc}
\frac{10}{12}&\frac{1}{12}&0&\cdots&0&0\\
\frac{1}{12}&\frac{10}{12}&\frac{1}{12}&\cdots&0&0\\
0&\frac{1}{12}&\frac{10}{12}&\cdots&0&0\\
\vdots&\vdots&\vdots&\ddots&\vdots&\vdots\\
0&0&0&\cdots&\frac{10}{12}&\frac{1}{12}\\
0&0&0&\cdots&\frac{1}{12}&\frac{10}{12}\\
\end{array}\right)_{M-1}, \quad \mathbf{S}_x=\left(\begin{array}{cccccc}
2&-1&0&\cdots&0&0\\
-1&2&-1&\cdots&0&0\\
0&-1&2&\cdots&0&0\\
\vdots&\vdots&\vdots&\ddots&\vdots&\vdots\\
0&0&0&\cdots&2&-1\\
0&0&0&\cdots&-1&2\\
\end{array}\right)_{M-1},
\end{eqnarray*}
and
 \begin{eqnarray*}
 &&\mathbf{M}_t^{\alpha}=\left(\begin{array}{llllll}
a_0^{(\alpha)}&0&0&\cdots&0&0\\
a_1^{(\alpha)}-a_0^{(\alpha)}&a_0^{(\alpha)}&0&\cdots&0&0\\
a_2^{(\alpha)}-a_1^{(\alpha)}&a_1^{(\alpha)}-a_0^{(\alpha)}&a_0^{(\alpha)}&\cdots&0&0\\
\vdots&\vdots&\vdots&\ddots&\vdots&\vdots\\
a_{N-2}^{(\alpha)}-a_{N-3}^{(\alpha)}&a_{N-3}^{(\alpha)}-a_{N-4}^{(\alpha)}&a_{N-4}^{(\alpha)}-a_{N-5}^{(\alpha)}&\cdots&a_0^{(\alpha)}&0\\
a_{N-1}^{(\alpha)}-a_{N-2}^{(\alpha)}&a_{N-2}^{(\alpha)}-a_{N-3}^{(\alpha)}&a_{N-3}^{(\alpha)}-a_{N-4}^{(\alpha)}&\cdots&a_1^{(\alpha)}-a_0^{(\alpha)}&a_0^{(\alpha)}\\
\end{array}\right)_N,\\
&&\mathbf{M}_t^{\beta}=\left(\begin{array}{llllll}
2b_0^{(\beta)}&0&0&\cdots&0&0\\
2b_1^{(\beta)}-2b_0^{(\beta)}&b_0^{(\beta)}&0&\cdots&0&0\\
2b_2^{(\beta)}-2b_1^{(\beta)}+b_0^{(\beta)}&b_1^{(\beta)}-2b_0^{(\beta)}&b_0^{(\beta)}&\cdots&0&0\\
2b_3^{(\beta)}-2b_2^{(\beta)}+b_1^{(\beta)}&b_2^{(\beta)}-2b_1^{(\beta)}+b_0^{(\beta)}&b_1^{(\beta)}-2b_0^{(\beta)}&\cdots&0&0\\
\vdots&\vdots&\vdots&\ddots&\vdots&\vdots\\
2b_{N-2}^{(\beta)}-2b_{N-3}^{(\beta)}+b_{N-4}^{(\beta)}&b_{N-3}^{(\beta)}-2b_{N-4}^{(\beta)}+b_{N-5}^{(\beta)} &b_{N-4}^{(\beta)}-2b_{N-5}^{(\beta)}+b_{N-6}^{(\beta)} &\cdots&b_0^{(\beta)}&0\\
2b_{N-1}^{(\beta)}-2b_{N-2}^{(\beta)}+b_{N-3}^{(\beta)}&b_{N-2}^{(\beta)}-2b_{N-3}^{(\beta)}+b_{N-4}^{(\beta)} &b_{N-3}^{(\beta)}-2b_{N-4}^{(\beta)}+b_{N-5}^{(\beta)} &\cdots&b_1^{\beta}-2b_0^{(\beta)}&b_0^{(\beta)}\\
\end{array}\right)_N,\\
\end{eqnarray*}

Thus, the difference scheme \eqref{b9}-\eqref{b11} can be equivalently reformulated as the following matrix equation:
\begin{eqnarray}\label{HL-matrix-eq}
(K_1\tau^{\beta-\alpha} \mathbf{M}_t^{\alpha}+K_2 \mathbf{M}_t^{\beta})\mathbf{u}\mathbf{M}_x+\frac{\tau^{\beta}}{h^2}\mathbf{u}\mathbf{S}_x=\mathbf{b} \mathbf{M}_x
\end{eqnarray}
where $(\mathbf{u})_{i,j}=u_{j}^{i},$ $i=1,\ldots N,$ $j=1,\ldots M-1,$  and the right hand side matrix $\mathbf{b}\in \mathbb{R}^{N\times (M-1)}$ is given by
\begin{eqnarray}
\mathbf{b}=\tau^{\beta}\mathbf{F} +2K_2\tau\left(\begin{array}{l}
b_0^{(\beta)} \\
b_1^{(\beta)} \\
\vdots\\
b_{N-2}^{(\beta)} \\
b_{N-1}^{(\beta)} \\
\end{array}\right)(\phi_1^1,\phi_2^1,\cdots,\phi^1_{M-2},\phi^1_{M-1} )
\end{eqnarray}
with $(\mathbf{F})_{i,j}=F_{j}^{i},$ $i=1,\ldots N,~j=1,\ldots, M-1.$
From Lemma \ref{lem-trimax-1}, we know that there exists a normalized  orthogonal matrix, $\mathbf{Q}_x\in \mathbb{R}^{(M-1)\times (M-1)},$  such that
\begin{eqnarray*}
&&\mathbf{\Lambda}_x^{(1)}=diag(\lambda_1^{(1)},\lambda_2^{(1)},\ldots, \lambda_{M-1}^{(1)}),\\
&&\mathbf{\Lambda}_x^{(2)}=diag(\lambda_1^{(2)},\lambda_2^{(2)},\ldots, \lambda_{M-1}^{(2)}),\\
&& \mathbf{M}_x\mathbf{Q}_x=\mathbf{Q}_x\mathbf{\Lambda}^{(1)}_x, \quad \mathbf{S_x}\mathbf{Q}_x=\mathbf{Q}_x\mathbf{\Lambda}_x^{(2)},
\end{eqnarray*}
More precisely, we have explicit representation
$$(\mathbf{Q}_x)_{i,j}= \sqrt{\frac{2}{M}}\sin( \frac{ij\pi}{M}),i,j=1,2,\ldots,M-1,$$
and
$$ \lambda_i^{(1)}=\frac{5}{6}+\frac{1}{6}\cos(\frac{i\pi}{M}),\quad  \lambda_i^{(2)}=2-2\cos(\frac{i\pi}{M}), \quad i= 1,2,\ldots,M-1.  $$
Hence, multiplying Equation \eqref{HL-matrix-eq} by $\mathbf{Q}_x$  gives
\begin{eqnarray}\label{HL-matrix-eq-2}
(K_1\tau^{\beta-\alpha} \mathbf{M}_t^{\alpha}+K_2 \mathbf{M}_t^{\beta})\mathbf{u}\mathbf{Q}_x\mathbf{\Lambda}_x^{(1)}+\frac{\tau^{\beta}}{h^2}\mathbf{u}\mathbf{Q}_x\mathbf{\Lambda}_x^{(2)}=\mathbf{b} \mathbf{Q}_x\mathbf{\Lambda}_x^{(1)}.
\end{eqnarray}
Let $\mathbf{v}\in \mathbb{R}^{N\times {M-1}}$ such that
$$\mathbf{u}=\mathbf{v}\mathbf{Q}_x. $$
Then, Equation \eqref{HL-matrix-eq-2} is equivalent to
\begin{eqnarray}\label{HL-matrix-eq-3}
(K_1\tau^{\beta-\alpha} \mathbf{M}_t^{\alpha}+K_2 \mathbf{M}_t^{\beta})\mathbf{v}\mathbf{\Lambda}_x^{(1)}+\frac{\tau^{\beta}}{h^2}\mathbf{v}\mathbf{\Lambda}_x^{(2)}=\mathbf{b} \mathbf{Q}_x\mathbf{\Lambda}_x^{(1)}=\mathbf{G}\mathbf{\Lambda}_x^{(1)},
\end{eqnarray}
where $ \mathbf{G}=\mathbf{b} \mathbf{Q}_x.$

Let $\mathbf{I}\in \mathbb{R}^{N\times N}$ be an identity matrix, and $\mathbf{v}_i=(v_{i1},v_{i2},\ldots,v_{iN} )^{T} $ and  $\mathbf{G}_i=(G_{i1},G_{i2},\ldots,G_{iN} )^{T}.$ Then the $i$-th column of the \eqref{HL-matrix-eq-3} can be written as
\begin{eqnarray}\label{HL-matrix-eq-4}
\biggl(K_1\tau^{\beta-\alpha} \mathbf{M}_t^{\alpha}+K_2 \mathbf{M}_t^{\beta}+\frac{\lambda_i^{(2)}}{\lambda_i^{(1)} }\frac{\tau^{\beta}}{h^2}\mathbf{I}\biggl)\mathbf{v}_i=\mathbf{G}_i,\quad i=1,2,\ldots,M-1,
\end{eqnarray}
which are equivalent to $M-1$ systems of the following form
\begin{eqnarray}\label{HL-matrix-eq-5}
 \left(\begin{array}{llllll}
c_0&0&0&\cdots&0&0\\
c_1&d_1&0&\cdots&0&0\\
c_2&d_2&d_1&\cdots&0&0\\
\vdots&\vdots&\vdots&\ddots&\vdots&\vdots\\
c_{N-2}&d_{N-2}&d_{N-3}&\cdots&d_1&0\\
c_{N-1}&d_{N-1}&d_{N-2}&\cdots&d_2&d_1\\
\end{array}\right)
 \left(\begin{array}{l}
e_0\\
e_1\\
e_2\\
\vdots\\
e_{N-2}\\
e_{N-1}\\
\end{array}\right)= \left(\begin{array}{l}
g_0\\
g_1\\
g_2\\
\vdots\\
g_{N-2}\\
g_{N-1}\\
\end{array}\right).
\end{eqnarray}

To fully employ the Toeplitz structure, we recast the above equation as
\begin{eqnarray}\label{HL-matrix-eq-5}
 \left(\begin{array}{cccccc}
d_1&0&\cdots&0&0\\
d_2&d_1&\cdots&0&0\\
\vdots&\vdots&\vdots&\ddots&\vdots&\vdots\\
d_{N-2}&d_{N-3}&\cdots&d_1&0\\
d_{N-1}&d_{N-2}&\cdots&d_2&d_1\\
\end{array}\right)
 \left(\begin{array}{c}
e_1\\
e_2\\
\vdots\\
e_{N-2}\\
e_{N-1}\\
\end{array}\right)= \left(\begin{array}{c}
g_1-c_1e^0\\
g_2-c_2e^0\\
\vdots\\
g_{N-2}-c_{N-2}e^0\\
g_{N-1}-c_{N-1}e^0\\
\end{array}\right)
\end{eqnarray}
and $e_0=f_0/c_0.$
The resulting lower triangular Toeplitz matrix equation \eqref{HL-matrix-eq-5} can be solved by the divide and conquer strategy proposed in \cite{Commenges1984}, where the authors gave the fast inversion of the lower triangular Toeplitz matrices. The idea of the divide and conquer strategy is also employed by the very recent paper \cite{Ke2015} for solving block triangular Toeplitz-like with tri-diagonal block systems from time-fractional partial differential equations. For the completeness of this paper, we sketch it below.

Suppose $Ax=b,$ where $A\in \mathbb{R}^{N}$ is a triangular Toeplitz matrix. For simplicity, we assume $N=2^{n}.$  Obviously, $A,$ $x$ and $b$ can be partitioned as follows:
\begin{equation}
 \left(\begin{array}{cc}
 A^{(1)}& 0\\
 C^{(1)} & A^{(1)}\\
 \end{array}\right)
  \left(\begin{array}{c}
  x^{(1)}\\
  x^{(2)}
  \end{array}\right)=
   \left(\begin{array}{c}
     b^{(1)}\\
  b^{(2)}
  \end{array}\right).
\end{equation}
Thus the original linear system can be equivalently transformed into two half size linear systems \begin{equation}
 \left\{\begin{array}{l}
 A^{(1)}x^{(1)}= b^{(1)}\\
 A^{(1)}x^{(2)}=b^{(2)}-C^{(1)}x^{(1)}\\
 \end{array}\right. .
\end{equation}
We note that $A^{(1)}$ is a still triangular Toeplitz matrix and $C^{(1)}$ a Toeplitz matrix. The matrix-vector product $C^{(1)}x^{(1)}$ can be computed efficiently by fast Fourier transform (FFT). The same procedure can be applied to solve both linear systems recursively.  Suppose that $\Theta_N$ is the number of  operations required to solve lower tri-diagonal linear system. The computation cost can be estimated below
$$ \Theta_N=2\Theta_{N/2}+\frac{N}{2}\log(\frac{N}{2}).$$
By this formula, we can derive the total operations in time direction as $\Theta_N=O(N\log^2N),$ which has great advantage than the  forward substitution  method with operations $O(N^2).$

 It should be pointed out that FFT can also applied to matrix-matrix or  matrix-vector product in spatial direction due to the periodicity and symmetry of the orthogonal matrix $\mathbf{Q}_x,$ with  operations of $O(M\log M)$  in space.  From aforementioned  analysis, we can derive that  the total operations of our proposed scheme  in space and time is $O(MN\log M \log^2N),$ which enjoys linearithmic complexity.

  As to two dimensional case, we can also use partial diagonalization technique to reduce the resulting  two-dimensional difference scheme into $(M_1-1)\times(M_2-1)$ linear  systems as \eqref{HL-matrix-eq-5}, which can be solved efficiently by the fast direct solver via the divide and conquer idea. Indeed, it suffices to partially diagonalize in the $y$ direction or $x$ direction first, then all the subsequent procedure are as the same as the one-dimensional case.  We will not dwell on the fast solver anymore and the robustness and accuracy of the scheme \eqref{b9-2D}-\eqref{b11-2D}  will be illustrated by  numerical example 7.3.

\section{Numerical examples}
\setcounter{equation}{0}
 In order to illustrate the behavior of our proposed  numerical method and demonstrate the effectiveness of our theoretical analysis, several  examples are
presented.

\begin{example}
As the first example, consider the problem
\begin{equation}\label{HL-15}
\  _0^{C}D_t^{\alpha_1}u(x,t)+\  _0^{C}D_t^{\alpha_2}u(x,t)+ u(x,t)= \partial_x^{2}u(x,t) +f(x,t),\quad 0< x < \pi,~0<t\leq 1,
\end{equation}
with  the source term and initial-boundary value conditions satisfying that  the equation admits the solution
$$u(x,t)=\sin (x) (t^3+t+1). $$
\end{example}
 Convergence orders, errors and CPU time of the compact scheme \eqref{b9}-\eqref{b11} are  examined.
Let $$E(h,\tau)=\max_{1\leq i\leq M-1}|u(x_i,t_N)-u_i^N|,$$
where $u(x_i,t_N) $ represents the exact solution and $u^N$ is the  numerical solution with the step sizes  $h$ and $\tau$ at $t_N=1.$
In simulation, we take parameter $\alpha_2=\alpha_1+1,$  under which Equation \eqref{HL-15} becomes the well-known time fractional telegraph equation. In this example, we use the proposed scheme \eqref{b9}-\eqref{b11} to solve this problem and implement it by the fast solver presented in Section 6.  The convergence order is computed by $$ \mbox{Order} = \frac{\log(E(\tau,h_1)/E(\tau,h_2) )}{\log(h_1/h_2)}  $$

From Tables \ref{table1}-\ref{table2}, it is clear to see that the convergence order is one in temporal direction and four in spatial direction, respectively,  which is in good agreement with the theoretical analysis. In addition, we can also see from the columns of 'CPU(s)' that the proposed  fast solver has a linearithmic complexity, which is in agreement with our prediction.

\begin{example}
Next, consider the following distributed order time-fractional diffusion-wave equation:
\begin{eqnarray*}
&&\int_{0}^2\Gamma(4-\alpha)\ _0^{C}D_t^{\alpha}u(x,t)d\alpha=\partial_x^2u(x,t)+\sin (x) (t^3+\frac{6t^3-6t}{\log t}),\quad 0<x<\pi,\quad 0<t\leq 1,
\end{eqnarray*}
with the exact solution  $u(x,t)=\sin(x) t^3.$
\end{example}
Setting
$$ \tilde{E}(h,\tau,\sigma)=\max_{1\leq i\leq M-1}|u(x_i,t_N,\sigma)-u_i^N|,$$ the errors and convergence orders are displayed in Table \ref{table3}. From the table, we can clearly see that the convergence orders are of first-order in time and fourth-order in space, respectively, and second-order with respect to the variable $\alpha,$  which again verifies the correctness of our theoretical results.

\begin{example}
Finally, consider the  two-dimensional problem \eqref{multi-dw-2D}-\eqref{boundary-2D} with exact solution
$$u(x,y,t)=\sin(x) \sin(y)(t^3+t+1)$$  on the domain $\Omega =(0,\pi)\times (0,\pi).$
\end{example}

\begin{figure}
\centering
\includegraphics[width=1\textwidth]{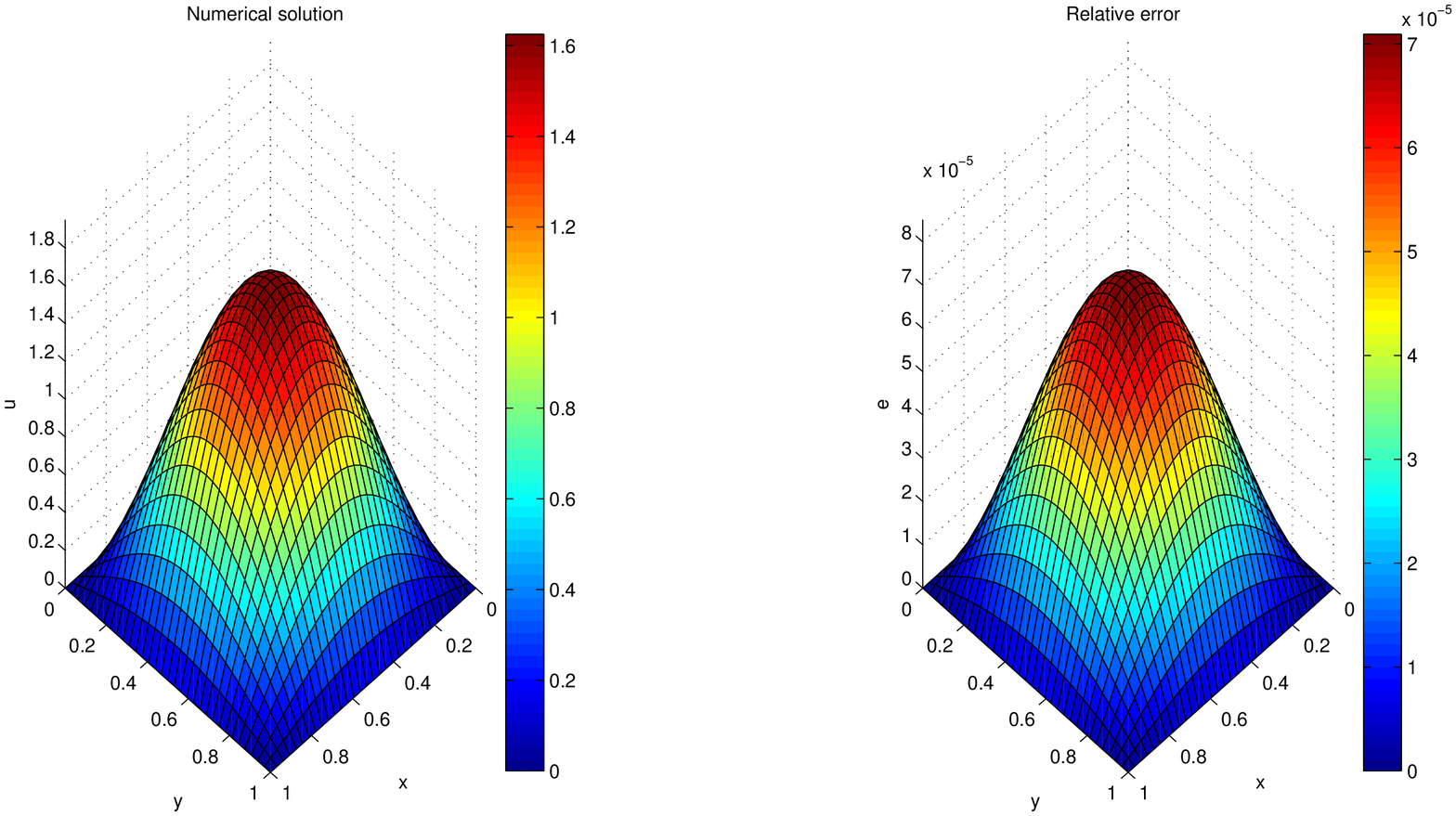}
\caption{Numerical solution and relative error  at $T=0.5$ with $\alpha_1=0.55$ and $\alpha_2=1.1$ (Ex. 8.3)} \label{fig:1}
\end{figure}
\begin{figure}
\centering
\includegraphics[width=1\textwidth]{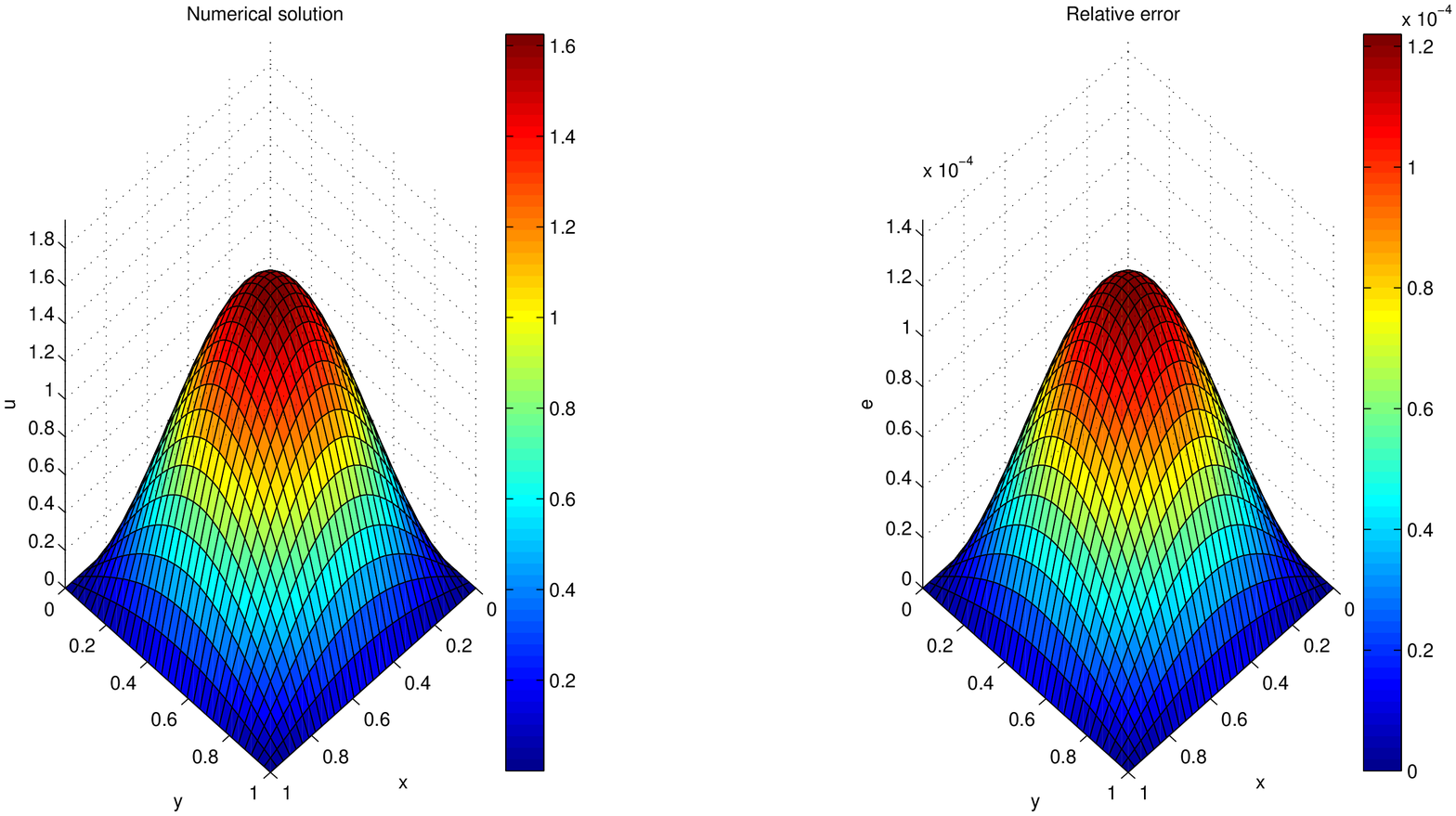}
\caption{Numerical solution and relative error  at $T=0.5$ with $\alpha_1=0.75$ and $\alpha_2=1.5$ (Ex. 8.3)} \label{fig:2}
\end{figure}
\begin{figure}
\centering
\includegraphics[width=1\textwidth]{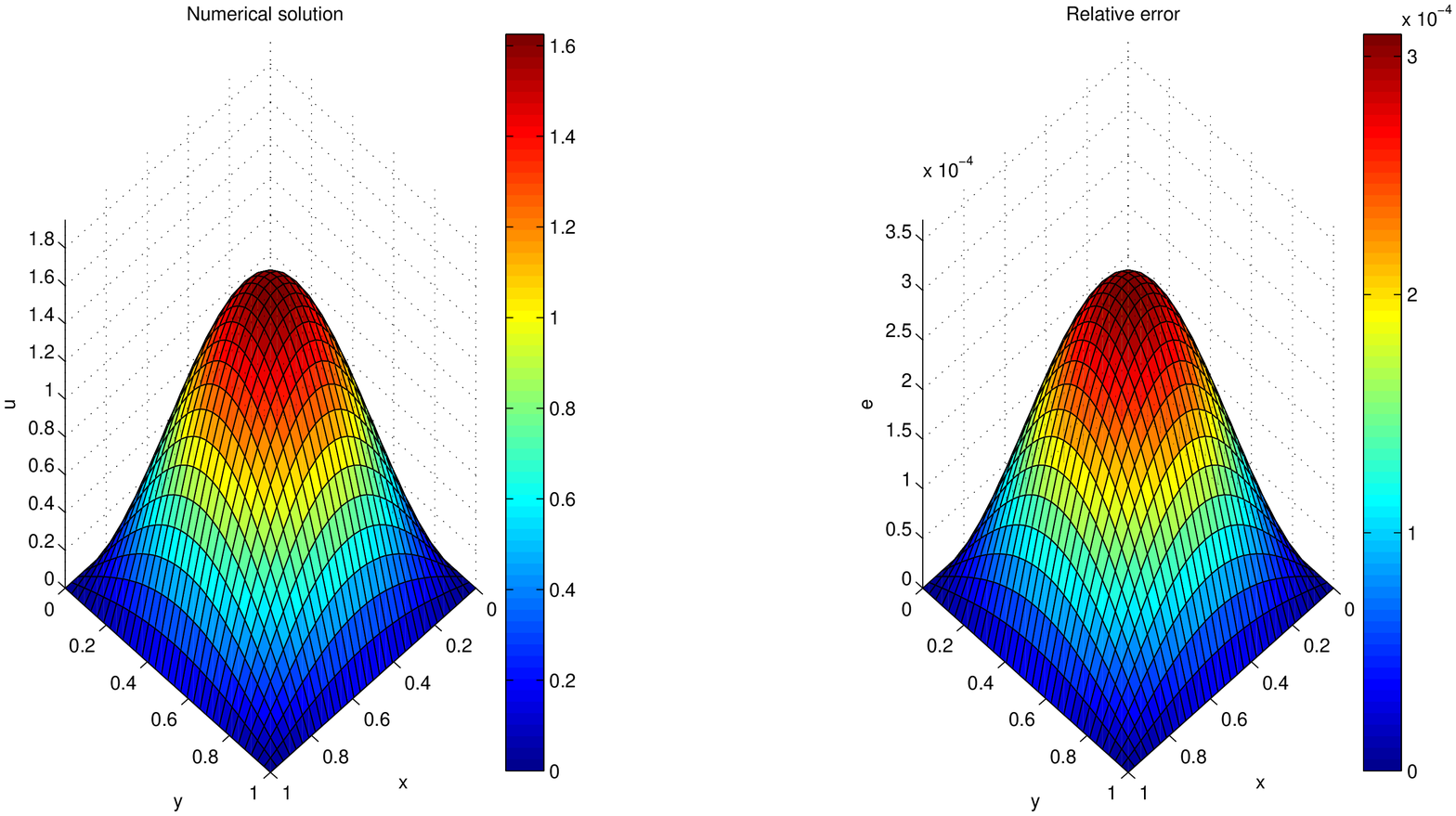}
\caption{Numerical solution and relative error  at $T=0.5$ with $\alpha_1=0.95$ and $\alpha_2=1.9$ (Ex. 8.3)} \label{fig:3}
\end{figure}

Take  $T=0.5.$  We take $h=1/32$ and $\tau=1/1024.$  We compute this problem by the scheme \eqref{b9-2D}-\eqref{b11-2D}.  Similar to  one-dimensional case,   we first  employ partial diagonalization technique to deal with $y$ direction, then we can compute the solution   following the exactly same procedure as one dimensional problem.
From Figure \ref{fig:1}-\ref{fig:3}, we can see that the numerical solution is perfectly matched with  the exact solution when fractional orders take different values.  This further verifies the robustness of the scheme  \eqref{b9-2D}-\eqref{b11-2D} and  the correctness of our theoretical results.

\section{Conclusion, remarks and discussion}
Numerical methods for multi-term  time-fractional  diffusion-wave equations are not fully developed yet up to now. The main goal of our work is to investigate the numerical solution and provide the complete error analysis to a class of mixed  diffusion-wave equations on bounded domain.  First of all,    a first-order approximation to the Caputo fractional derivative with order belonging to $(1,2)$ by modifying the $L2$ approximation is strictly  derived in this paper. Then the new approximation, combining with the $L1$ approximation to discretize the time derivative,  has been applied to solve the  multi-term time-fractional mixed diffusion-wave equations.    For the error analysis,  the discrete version of   the  fractional Poncar$\acute{e}$-Friedrichs Sobolev embedding   inequality  is provided and proved strictly, which plays an essential role in the numerical analysis. Then,    unconditional convergence of the compact difference scheme has been proved by the energy method, and stability can be obtained immediately with the consistency of the scheme.  Next, the distribution order equations and two-dimensional  extensions are considered.   Particularly,  we investigated    the distributed order fractional equations with orders belonging to (0,2). In addition, a practical  fast solver   with  linearithmic complexity are presented.   Finally, several numerical examples have been given to show the effectiveness and correctness of our proposed schemes.

Finite difference method is widely applied to solve time-fractional  equations  mainly because it is  easy to implement,  and fast solver can be developed without too much difficulty,  owning  to the Toeplitz structure of the resulting linear system.   Besides the popular Lagrange interpolant method, $L1$ and $L2$ approximation, there are other two commonly used ones, Gr\"{u}nwald-Letnikov approximation and  Lubich's fractional liner multi-step methods.  Gr\"{u}nwald-Letnikov approximation
\begin{equation}\label{a2}
\sum_{k=0}^{[t/\tau]}g^{(\alpha)}_ky(t-k\tau)= {_0^{RL}}D_t^{\alpha}y(t)+\mathcal{O}(\tau),\quad g^{(\alpha)}_k=(-1)^k\binom{\alpha}{k},
\end{equation}
 as an another way to discretize fractional differential and integral operators,  is based on the straightforward generalization of concepts from classical calculus to the fractional ones.  Since the coefficients $g^{(\alpha)}_k$ have the same properties as $a^{(\alpha)}_k$ in $L1$ approximation, which can be validated readily,  the interested  readers can also apply Gr\"{u}nwald-Letnikov formula to time-fractional mixed diffusion-wave equations and employ a similar analysis following the proof of this work. Lubich's fractional linear multi-step methods have fast convergence compared to two former  methods but  require    high regularity of the solution and correction of the error,   which limit their applications in practice.

There are also other directions deserved  future research. In this paper, we only consider the equations with sufficiently smooth solution,  $\partial _t^2u(\cdot,t)\in C^2[0,T].$  That is, the second order derivatives can be bounded.  We have not touched upon the low regularity case.  For example,  consider the problem \eqref{HL-15} with low regularity solution, $u(x,t)=\sin (x) t^\nu. $ From Table \ref{table5}, we can clearly see that convergence rate deteriorates with decrease of regularity of the solution and it behaves like $O(\tau^{\nu-1})$ in time. Actually,  this  can be derived directly from the approximation formula \eqref{denotion-be}.
  Due to the potential  singularity of the solution, nonuniform mesh based finite difference methods   have gained considerable attention for the time-fractional sub-diffusion     recently.  However, the corresponding  results about time-fractional mixed diffusion-wave equation is still scarce.    And we will consider this problem and present an effective fast solver    in our future  work.

\section*{Acknowledgement}
The authors appreciate Prof. Zhi-zhong Sun for providing us useful suggestions through careful reading.
The research   is
supported by National Natural Science Foundation of China (No.
11271068) and by the Fundamental Research Funds for the Central Universities and the Research and Innovation Project for College Graduates of Jiangsu Province (Grant No.: KYLX\_0081).  The research  is also partially  supported by the Scientific Research Foundation of Graduate School of Southeast University.
G. Lin would like to thank the support by NSF Grant DMS-1115887, and the U.S. Department of Energy, Office of Science, Office of Advanced Scientific Computing Research, Applied Mathematics program as part of the Collaboratory on Mathematics for Mesoscopic Modeling of Materials, and Multifaceted Mathematics for Complex Energy Systems.

\begin{table}[!htb]
\centering
\caption{ Temporal convergence orders, errors and CPU time  of the scheme \eqref{b9}-\eqref{b11} with fixed stepsize $h=1/16$ (Ex. 8.1) }\label{table1}
\vspace{0.1in}
\newsavebox{\tablebox}
\begin{lrbox}{\tablebox}
\begin{tabular}{cccccccccc}
\hline\hline
  & \multicolumn{3}{c}{ $\alpha_1=0.2,~\alpha_2=1.2$}&   \multicolumn{3}{c}{ $\alpha_1=0.5,~\alpha_2=1.5$}  &\multicolumn{3}{c}{ $\alpha_1=0.7,~\alpha_2=1.7$} \\
\cline{2-4} \cline{5-7} \cline{8-10}
    $N$ &  $E(h,\tau) $ & $Order$ & CPU(s) & $E(h,\tau) $ & $Order$ & CPU(s) &$E(h,\tau) $ & $Order$ & CPU(s) \\
     \hline
     16     & 4.4722e-2     &    -- &  0.0411   &  5.7604e-2     &    --  &0.0216    & 7.5149e-2       &    -- &   0.0199    \\
      32& 2.2796e-2     &   0.9722& 0.0694    & 2.8439e-2     &   1.0183& 0.0402    & 3.6662e-2     &   1.0355  &   0.0382     \\
   64  & 1.1489e-2     &   0.9885 &0.0802     &1.3953e-2     &  1.0273 &0.0810    & 1.7712e-2     &   1.0496 &    0.0789     \\
     128& 5.7626e-3     &   0.9955& 0.1596  & 6.8480e-3     &   1.0268  &  0.1555& 8.5411e-3     &   1.0522 & 0.1816 \\
\hline
\hline
\end{tabular}
\end{lrbox}
\scalebox{1.0}{\usebox{\tablebox}} 
\end{table}

\begin{table}[!htb]
\centering
\caption{ Spatial convergence orders, errors and CPU time the scheme \eqref{b9}-\eqref{b11} with fixed stepsize $\tau=1/2^{20}$ (Ex. 8.1) }\label{table2}
\vspace{0.1in}
\newsavebox{\tableboxpeng}
\begin{lrbox}{\tableboxpeng}
\begin{tabular}{cccccccccc}
\hline\hline
  & \multicolumn{3}{c}{ $\alpha_1=0.2,~\alpha_2=1.2$}&   \multicolumn{3}{c}{ $\alpha_1=0.5,~\alpha_2=1.5$}  &\multicolumn{3}{c}{ $\alpha_1=0.7,~\alpha_2=1.7$} \\
\cline{2-4} \cline{5-7} \cline{8-10}
    $M$ &  $E(h,\tau) $ & $Order$ & CPU(s) & $E(h,\tau) $ & $Order$ & CPU(s) &$E(h,\tau) $ & $Order$ & CPU(s) \\
     \hline
     4     & 1.0743e-3    &    --  &  378.44 &1.0073e-3     &    --  &  365.83  & 9.2285e-4        &    --   & 362.97   \\
      6& 2.1008e-4     &   4.0248 &  564.18  &1.9709e-4     &   4.0234&  544.93 & 1.7724e-4     &   4.0694   &   528.73  \\
   8  & 6.6644e-5     &   3.9910  & 701.06  &6.2494e-5     &   3.9927 & 714.85 & 5.3393e-5     &   4.1706 &720.86\\
     10& 2.7652e-5     &   3.9421  &872.26& 2.5945e-5     &   3.9396 &  905.42 & 2.0317e-5     &   4.3301 &868.36\\
\hline
\hline
\end{tabular}
\end{lrbox}
\scalebox{1.0}{\usebox{\tableboxpeng}} 
\end{table}

\begin{table}[!htb]
\centering
\caption{ Convergence orders and errors of the scheme \eqref{HL-12}-\eqref{HL-14} (Ex. 8.2) }\label{table3}
\vskip0.1in
\tabcolsep0.16in
\begin{tabular}{ccccccccc}
\hline\hline
   \multicolumn{3}{c}{ $h=1/16$, $\sigma=1/16$}&   \multicolumn{3}{c}{ $\tau=1/2^{20}$, $\sigma=1/128$}  &\multicolumn{3}{c}{ $\tau=1/2^{16}$, $h=1/16$ } \\
\cline{1-3} \cline{4-6} \cline{7-9}
    $N$& $\tilde{E}(h,\tau,\sigma) $ & $Order$ &M & $\tilde{E}(h,\tau,\sigma) $ & $Order$ &J &$\tilde{E}(h,\tau,\sigma) $ & $Order$  \\
      \hline
    16     & 4.4722e-2     &    --  &4&8.5302e-5  &  --& 2     & 2.8526e-3     &    -- \\
      32& 2.2796e-2     &   0.9722  &6&1.5113e-5  &  4.2683& 4& 7.3243e-4     &   1.9615 \\
   64  & 1.1489e-2     &   0.9885   &8& 3.6059e-6 & 4.9811&6  & 3.3198e-4     &     1.9516 \\
     128& 5.7626e-3     &   0.9955  &10& 8.4501e-7&  6.5024 &8& 1.9127e-4     &    1.9167 \\
    \hline
\hline
\end{tabular}
\end{table}

\begin{table}[!htb]
\centering
\caption{ Temporal convergence orders and errors of the scheme \eqref{b9}-\eqref{b11} with fixed stepsize $h=1/16$, }\label{table5}
\vspace{0.1in}
\newsavebox{\tableboxg}
\begin{lrbox}{\tableboxg}
\begin{tabular}{ccccccc}
\hline\hline
  & \multicolumn{2}{c}{ $\alpha_1=0.6,~\alpha_2=1.2,~\nu=1.2$}&   \multicolumn{2}{c}{ $\alpha_1=0.75,~\alpha_2=1.5,~\nu=1.5$}  &\multicolumn{2}{c}{ $\alpha_1=0.85,~\alpha_2=1.7,~\nu=1.7$} \\
\cline{2-3} \cline{4-5} \cline{6-7}
    $N$ &  $E(h,\tau) $ & $Order$ &  $E(h,\tau) $ & $Order$  &$E(h,\tau) $ & $Order$  \\
     \hline
  16       & 1.9716e-1     &     --  & 1.0473e-1     &  --& 6.6505e-2     &     --   \\
 32       & 1.7182e-1     &   0.1984 & 7.5392e-2     &   0.4742  & 4.3388e-2     &   0.6161 \\
 64       & 1.4973e-1     &   0.1985  & 5.3799e-2     &   0.4868 & 2.7687e-2     &   0.6481  \\
 128       & 1.3045e-1     &   0.1989   & 3.8218e-2     &   0.4933   & 1.7431e-2     &   0.6676  \\
\hline
\hline
\end{tabular}
\end{lrbox}
\scalebox{1.0}{\usebox{\tableboxg}} 
\end{table}


 \end{document}